\documentclass[10pt]{article}

\usepackage{amsmath}
\usepackage{amsthm}
\usepackage{amssymb}
\usepackage{amscd}
\usepackage{amsfonts}  
\usepackage{makeidx}
\usepackage{hyperref}
\usepackage[all]{xy} \SelectTips{cm}{}
\usepackage{makeidx}     
\usepackage{graphicx}    
\usepackage{multicol}    
\usepackage{pdfsync}

\DeclareMathAlphabet\EuR{U}{eur}{m}{n}
\SetMathAlphabet\EuR{bold}{U}{eur}{b}{n}

\makeindex             

\theoremstyle{plain}
\newtheorem{theorem}{Theorem}[section]
\newtheorem{lemma}[theorem]{Lemma}

\theoremstyle{definition}

\newtheorem{example}[theorem]{Example}

\newtheorem{remark}[theorem]{Remark}

{\catcode`@=11\global\let\c@equation=\c@theorem}

\newcommand{\comsquare}[8]                   
{\begin{CD}
#1 @>#2>> #3\\
@V{#4}VV @V{#5}VV\\
#6 @>#7>> #8
\end{CD}
}

\newcommand{\xycomsquare}[8]                   
{\xymatrix
{#1 \ar[r]^{#2} \ar[d]^{#4} &
#3 \ar[d]^{#5}  \\
#6\ar[r]^{#7} &
#8
}
}

\newcommand{\calfin}{\mathcal{FIN}}

\newcommand{\calvcyc}{\mathcal{VCY}}
\newcommand{\caltr}{\{ \! 1 \! \}}

\newcommand{\calall}{{\cala\!\call\!\call}}
\newcommand{\calmicyc}{\mathcal{MICY}}

\newcommand{\cala}{{\cal A}}

\newcommand{\calf}{{\cal F}}
\newcommand{\calg}{{\cal G}}
\newcommand{\calh}{{\cal H}}

\newcommand{\call}{{\cal L}}
\newcommand{\calm}{{\cal M}}

\newcommand{\cals}{{\cal S}}

\newcommand{\IQ}{{\mathbb Q}}
\newcommand{\IR}{{\mathbb R}}

\newcommand{\IZ}{{\mathbb Z}}

\newcommand{\bfK}{{\mathbf K}}

\newcommand{\bfL}{{\mathbf L}}

\newcommand{\curs}{\EuR}

\newcommand{\GROUPOIDS}{\curs{GROUPOIDS}}

\newcommand{\MODULES}{\curs{MODULES}}
\newcommand{\Or}{\curs{Or}}

\newcommand{\SPECTRA}{\curs{SPECTRA}}

\newcommand{\aut}{\operatorname{aut}}
\newcommand{\CAT}{\operatorname{CAT}}

\newcommand{\coker}{\operatorname{coker}}

\newcommand{\cyl}{\operatorname{cyl}}

\newcommand{\id}{\operatorname{id}}

\newcommand{\ind}{\operatorname{ind}}

\newcommand{\Nil}{\operatorname{Nil}}

\newcommand{\per}{\operatorname{per}}
\newcommand{\pr}{\operatorname{pr}}

\newcommand{\res}{\operatorname{res}}

\newcommand{\UNil}{\operatorname{UNil}}

\newcommand{\Wh}{\operatorname{Wh}}

\newcommand{\pt}{\{\bullet\}}
\newcommand{\NK}{N\!K}

\newcommand{\EGF}[2]{E_{#2}(#1)}                   
\newcommand{\eub}[1]{\underline{E}#1}              
\newcommand{\edub}[1]{\underline{\underline{E}}#1} 

\newcommand{\higherlim}[3]{{\setbox1=\hbox{\rm lim}
     \setbox2=\hbox to \wd1{\leftarrowfill} \ht2=0pt \dp2=-1pt
     \mathop{\vtop{\baselineskip=5pt\box1\box2}}
     _{#1}}^{#2}#3}

\newcommand{\version}[1]                       
{\begin{center} last edited on #1\\
last compiled on \today\\
name of tex-file: \jobname
\end{center}
}

\newcounter{commentcounter}

\begin{document}

\typeout{----------------------------  klva_final.tex  ----------------------------}

\title{On the $K$- and $L$-theory of hyperbolic and  virtually finitely generated abelian groups}
\author{Wolfgang L\"uck\thanks{\noindent email:
wolfgang.lueck@him.uni-bonn.de} \protect\\
http://www.him.uni-bonn.de/lueck\protect\\
Mathematisches Institut\\ Universit\"at Bonn\\
Endenicher Allee~60\\ 53115 Bonn\\Germany
\and
David Rosenthal\thanks{\noindent email: rosenthd@stjohns.edu}\\
Department of Mathematics and Computer Science\\
St. John's University\\
8000 Utopia Pkwy\\
Jamaica, NY 11439\\
USA\\
}
\date{May 2012}
\maketitle


\typeout{-----------------------  Abstract  ------------------------}

\begin{abstract}
We investigate the algebraic $K$- and $L$-theory of the group ring
$RG$, where $G$ is a hyperbolic or virtually finitely generated
abelian group and $R$ is an associative ring with unit.
\\[2mm]
Key words: Algebraic $K$- and $L$-theory, hyperbolic groups,
virtually finitely generated abelian groups.\\
Mathematics Subject Classification 2000: 19D99, 19G24, 19A31, 19B28.
\end{abstract}


\typeout{----------------------------------------------   Introduction ----------------------------------}
\section*{Introduction}

The goal of this paper is to compute the algebraic $K$- and $L$-groups
of group rings $RG$, where $G$ is a hyperbolic group or a virtually
finitely generated abelian group and $R$ is an associative ring with
unit (and involution).

The Farrell-Jones Conjecture for a group $G$ and a ring $R$
predicts that the assembly maps
\begin{eqnarray*}
H_n^G\bigl(\edub{G};\bfK_R\bigr) 
& \xrightarrow{\cong} & 
H_n^G\bigl(G/G;\bfK_R\bigr) = K_n(RG);
\\
H_n^G\bigl(\edub{G};\bfL_R^{\langle - \infty \rangle}\bigr) 
& \xrightarrow{\cong} & 
H_n^G\bigl(G/G;\bfL_R^{\langle - \infty \rangle}\bigr) = L_n^{\langle - \infty \rangle}(RG),
\end{eqnarray*}
induced by the projection $\edub{G} \to G/G$, are bijective for every integer $n$. This conjecture, introduced by Farrell and Jones in their groundbreaking
paper~\cite{Farrell-Jones(1993a)}, has many consequences. Knowing that it is
true for a given group implies several other well-known conjectures for that
group, such as the ones due to Bass, Borel, Kadison and Novikov. The Farrell-Jones
Conjecture also helps to calculate the $K$- and $L$-theory of group rings, since
homology groups are equipped with tools such as spectral sequences that can
simplify computations.

Recently it has been established that the Farrell-Jones Conjecture is true for
word hyperbolic groups and virtually finitely generated abelian groups. Using
this fact, we are able to compute the $K$- and $L$-theory of their group rings
in several cases by analyzing the left-hand side of the assembly map. The key
ingredients used to compute these groups are the \emph{induction structure} that
equivariant homology theories possess~\cite[Section~1]{Lueck(2002b)} and the
work of L\"{u}ck-Weiermann~\cite{Lueck-Weiermann(2012)}, which investigates when a universal space for a given group $G$ and a given family of subgroups
$\calf$ can be constructed from a universal space for $G$ and a smaller family
$\calf'\subseteq \calf$.

Even in basic situations determining the $K$- and $L$-groups is
difficult. However, we are able to handle hyperbolic
groups. The favorite situation is when the hyperbolic group $G$ is torsion-free and $R = \IZ$, in which case $K_n(\IZ G) = 0$ for $n \le -1$, the reduced projective class group $\widetilde{K}_0(\IZ G)$ and the Whitehead group $\Wh(G)$ vanish, and $K_n(\IZ G)$ is computed by $H_n(BG;\bfK(\IZ))$, i.e., the homology with
coefficients in the $K$-theory spectrum $\bfK(\IZ)$ of the classifying space
$BG$. Moreover, the $L$-theory $L_n^{\langle i \rangle}(\IZ G)$ is independent
of the decoration $\langle i \rangle$ and is given by $H_n(BG;\bfL(\IZ))$, where
$\bfL(\IZ)$ is the (periodic) $L$-theory spectrum of $\bfL(\IZ)$. Recall that
$\pi_n(\bfL(\IZ))$ is $\IZ$ if $n \equiv 0\mod 4$, $\IZ/2$ if if $n \equiv 2
\mod 4$, and vanishes otherwise. We also give a complete answer for $G = \IZ^d$. However, for a group $G$
appearing in an exact sequence $1 \to \IZ^d \to G \to Q \to 1$ for a finite
group $Q$, our computations can only be carried out under the additional assumption that the conjugation action of $Q$ on $\IZ^d$ is free away from $0$ or
that $Q$ is cyclic of prime order.

In Section~\ref{sec:statements} precise statements of our results are given, as
well as several examples. Section~\ref{sec:Background} contains the necessary
background for the proofs, which are presented in Section~\ref{sec:proofs}.

The paper was supported by the Sonder\-forschungs\-be\-reich SFB 878 -- Groups,
Geometry and Actions --, the Leibniz-Preis of the first author, and the
Fulbright Scholars award of the second author. The authors thank the referee for several 
useful comments.

\tableofcontents


\typeout{----------------------------   Section 1: Statement of Results --------------------------------------}

\section{Statement of Results}\label{sec:statements}
Let $K_n(RG)$ denote the \emph{algebraic $K$-groups} of the group ring $RG$ in
the sense of Quillen for $n \ge 0$ and in the sense of Bass for $n \le -1$, and
let $L_n^{\langle - \infty \rangle}(RG)$ be the \emph{ultimate lower quadratic
$L$-groups} of $RG$ (see Ranicki~\cite{Ranicki(1992a)}). When considering
$L$-theory, we will tacitly assume that $R$ is a ring with an involution. Let
$\NK_n(R)$ denote the \emph{Bass-Nil-groups} of $R$, which are defined as the cokernel of the map
$K_n(R) \to K_n(R[x])$.  Recall that the
Bass-Heller-Swan decomposition says
\begin{eqnarray}
K_n(R\IZ) & \cong & K_n(R) \oplus K_{n-1}(R) \oplus \NK_n(R) \oplus \NK_n(R).
\label{Bass-Heller-Swan}
\end{eqnarray}
If $R$ is a regular ring, then $\NK_n(R) = 0$ for every $n \in \IZ$ (see for
instance Rosenberg~\cite[Theorems~3.3.3 and 5.3.30]{Rosenberg(1994)}).

The \emph{$n$-th Whitehead group} $\Wh_n(G;R)$ is defined as
\begin{eqnarray*}
\Wh_n(G;R) & := & 
H_n^G\bigl(EG \to \pt;\bfK_R\bigr),
\end{eqnarray*}
where $H_n^G\bigl(EG \to \pt;\bfK_R\bigr)$ is the relative term in the
equivariant homology theory $H_*^G\bigl(-;\bfK_R\bigr)$ associated to the
equivariant map $EG \to \pt$ (see~\ref{calh_nG(f)}).  Thus, $\Wh_n(G;R)$ fits
into a long exact sequence
\begin{multline*}
\cdots \to H_n\bigl(BG;\bfK(R)\bigr) \to K_n(RG) \to \Wh_n(G;R)
\\
\to H_{n-1}\bigl(BG;\bfK(R)\bigr) \to K_{n-1}(RG) \to \cdots,
\end{multline*}
where $\bfK(R)$ is the non-connective $K$-theory spectrum associated to $R$ and
$H_*\bigl(-;\bfK(R)\bigr)$ is the associated homology theory. Notice that
$\Wh_1(G;\IZ)$ agrees with the classical Whitehead group $\Wh(G)$. Whitehead
groups arise naturally when studying $h$-cobordisms, pseudoisotopy, and
Waldhausen's A-theory. Their geometric significance can be reviewed, for
example, in Dwyer-Weiss-Williams~\cite[Section~9]{Dwyer-Weiss-Williams(2003)}
and L\"uck-Reich~\cite[Section~1.4.1]{Lueck-Reich(2005)}, where additional
references can also be found.  When $G=\IZ$, it follows
from~\eqref{Bass-Heller-Swan} and the fact that
$H_n^{\IZ}\bigl(E\IZ;\bfK_R\bigr)\cong K_n(R) \oplus K_{n-1}(R)$ that there is
an identification
\begin{eqnarray}
H_n^{\IZ}\bigl(E\IZ \to \pt;\bfK_R\bigr) & \cong & \NK_n(R) \oplus \NK_n(R).
\label{H_nZ(EZ_to_pt;bfK_R)_cong_NK_n(R)_oplus_NK_n(R)}
\end{eqnarray}

Define the \emph{periodic $n$-th structure group with decoration $\langle -\infty \rangle$} to be
\begin{eqnarray*}
\cals_n^{\per,\langle -\infty \rangle}(G;R) 
& := &
H_n^G\bigl(EG \to \pt;\bfL_R^{\langle -\infty\rangle}\bigr).
\end{eqnarray*}
These groups fit into the periodic version of the long exact
surgery sequence with decoration $\langle -\infty \rangle$,
\begin{multline*}
\cdots \to H_n\bigl(BG;\bfL^{\langle -\infty \rangle}(R)\bigr) \to
L_n^{\langle -\infty \rangle}(RG) \to \cals_n^{\per,\langle -\infty
 \rangle}(G;R)
\\
\to H_{n-1}\bigl(BG;\bfL^{\langle -\infty \rangle}(R)\bigr) \to
L_{n-1}^{\langle -\infty \rangle}(RG) \to \cdots,
\end{multline*}
where $\bfL^{\langle -\infty \rangle}(R)$ is the spectrum
whose homotopy groups are the ultimate lower quadratic $L$-groups.
This periodic surgery sequence (with a different decoration) for $R = \IZ$ appears 
in the classification of ANR-homology manifolds in 
Bryant-Ferry-Mio-Weinberger~\cite[Main Theorem]{Bryant-Ferry-Mio-Weinberger(1996)}.
It is related to the algebraic surgery exact sequence and thus to the
classical surgery sequence (see Ranicki~\cite[Section~18]{Ranicki(1992)}).
For $n \in \IZ$ define
\begin{eqnarray}
\UNil^{\langle - \infty \rangle}_n(D_{\infty};R) & := & 
H_n^{D_{\infty}}\bigl(\eub{D_{\infty}}  
\to  \pt;\bfL_R^{\langle - \infty \rangle}\bigr).
\label{Unil-infty(D_infty;R)}
\end{eqnarray}
These groups are related to Cappell's $\UNil$-groups, as explained
in Section~\ref{sec:Background}.

Now we are able to state our main results.


\subsection{Hyperbolic groups}
\label{subsec:Hyperbolic_groups}

\begin{theorem}[Hyperbolic groups]
\label{the:hyperbolic_groups}
Let $G$ be a hyperbolic group in the sense of Gromov~\cite{Gromov(1987)}, and
let $\calm$ be a complete system of representatives of the conjugacy classes
of maximal infinite virtually cyclic subgroups of $G$.
\begin{enumerate}
\item \label{the:hyperbolic_groups:K} For each $n \in \IZ$ there is an
  isomorphism
  \begin{eqnarray*}
    H_n^G\bigl(\eub{G};\bfK_R\bigr) \oplus \bigoplus_{V \in \calm} 
    H_n^{V}\bigl(\underline{E}V \to \pt;\bfK_R\bigr) 
    & \xrightarrow{\cong} & K_n(RG);
  \end{eqnarray*}
\item \label{the:hyperbolic_groups:L} For each $n \in \IZ$ there is an
  isomorphism
  \begin{eqnarray*}
    H_n^G\bigl(\eub{G};\bfL_R^{\langle -\infty \rangle}\bigr) \oplus \bigoplus_{V \in \calm} 
    H_n^{V}\bigl(\underline{E}V \to \pt;\bfL_R^{\langle -\infty \rangle}\bigr) 
    &\xrightarrow{\cong} &
    L_n^{\langle -\infty \rangle}(RG),
  \end{eqnarray*}
  provided that there exists $n_0 \le -2$ such that $K_n(RV) = 0$ holds for
  all $n \le n_0$ and all virtually cyclic subgroups $V \subseteq G$. (The
  latter condition is satisfied if $R = \IZ$ or if $R$ is regular with $\IQ
  \subseteq R$.)
\end{enumerate}
\end{theorem}

A good model for $\eub{G}$ is given by the Rips complex of $G$ (see
Meintrup-Schick~\cite{Meintrup-Schick(2002)}).  Tools for computing
$H_n^G\bigl(\eub{G};\bfK_R\bigr)$ are the equivariant version of the
Atiyah-Hirzebruch spectral sequence (see
Davis-L\"uck~\cite[Theorem~4.7]{Davis-Lueck(1998)}), the $p$-chain spectral
sequence (see Davis-L\"uck~\cite{Davis-Lueck(2003)}) and equivariant Chern
characters (see L\"uck~\cite{Lueck(2002b)}).  More information about the groups
$H_n^{V}\bigl(\underline{E}V \to \pt;\bfK_R\bigr)$ and
$H_n^{V}\bigl(\underline{E}V \to \pt;\bfL^{\langle -\infty \rangle}_R\bigr)$ is
given in Section~\ref{sec:Background}.


\subsection{Torsion-free hyperbolic groups}\label{subsec:torsion-free_hyperbolic_groups}

The situation simplifies if $G$ is assumed to be a torsion-free hyperbolic group.

\begin{theorem}[Torsion-free hyperbolic groups]
\label{the:torsionfree_hyperbolic_groups}
Let $G$ be a torsion-free hyperbolic group, and let $\calm$ be a complete
system of representatives of the conjugacy classes of maximal infinite cyclic
subgroups of $G$.
\begin{enumerate}
\item \label{the:torsionfree_hyperbolic_groups:K} For each $n \in \IZ$ there
  is an isomorphism
  \begin{eqnarray*}
    H_n\bigl(BG;\bfK(R)\bigr) \oplus \bigoplus_{V \in \calm}  \NK_n(R) \oplus  \NK_n(R) 
    & \xrightarrow{\cong} & 
    K_n(RG);
  \end{eqnarray*}

\item \label{the:torsionfree_hyperbolic_groups:L} For each $n \in \IZ$ there
  is an isomorphism
  \begin{eqnarray*}
    H_n\bigl(BG;\bfL^{\langle -\infty \rangle}(R)\bigr)  
    & \xrightarrow{\cong} &
    L_n^{\langle -\infty \rangle}(RG).
  \end{eqnarray*}
\end{enumerate}
\end{theorem}

In particular, it follows that for a torsion-free hyperbolic group $G$ and a
regular ring $R$,
$$K_n(RG) = \{0\} \quad \text{for} \; n \le 1$$
and the obvious map
$$K_0(R) \xrightarrow{\cong} K_0(RG)$$
is bijective. Moreover, $\widetilde{K}_0(\IZ G)$, $\Wh(G)$ and $K_n(\IZ G)$ for
$n \le -1$ all vanish.

\begin{example}[Finitely generated free groups]
\label{exa:Finitely_generated_free_groups}
Let $F_r$ be the finitely generated  free group $\ast_{i=1}^r \IZ$ of
rank $r$. Since $F_r$ acts freely on a tree it is hyperbolic. By
Theorem~\ref{the:torsionfree_hyperbolic_groups},
\[
K_n(RF_r) 
\; \cong \; 
K_n(R) \oplus K_{n-1}(R)^r  \oplus \bigoplus_{V \in \calm}  \big(\NK_n(R)\oplus \NK_n(R)\big)\] 
and
\[
L_n^{\langle - \infty \rangle}(RF_r) 
\; \cong \;
L^{\langle - \infty \rangle}_n(R) \oplus L_{n-1}^{\langle - \infty \rangle}(R)^r,\]
where $\calm$ is a complete system of representatives of 
the conjugacy classes of maximal infinite cyclic subgroups of $F_r$.
\end{example}


\subsection{Finitely generated free abelian groups}
\label{subsec:IZd}

Before tackling the virtually finitely generated abelian case, we first consider finitely generated free abelian groups. The $K$-theory case of Theorem~\ref{the:Zd} below is also proved in~\cite{Davis(2008nil)}.

\begin{theorem}[$\IZ^d$]\label{the:Zd}
Let $d \ge 1$ be an integer. Let $\calmicyc$ be the set of maximal
infinite cyclic subgroups of $\IZ^d$.  Then there are isomorphisms
\begin{eqnarray*} 
\Wh_n(\IZ^d;R) 
& \cong  & 
\bigoplus_{C \in \calmicyc} \; \bigoplus_{i = 0}^{d-1} 
\bigl(\NK_{n-i}(R) \oplus \NK_{n-i}(R)\bigr)^{\binom{d-1}{i}};
\\
K_n(R[\IZ^d]) 
& \cong & 
\left(\bigoplus_{i = 0}^{d} K_{n-i}(R)^{\binom{d}{i}} \right) \oplus \Wh_n(\IZ^d;R);
\\
L_n^{\langle -\infty \rangle}(R[\IZ^d]) 
& \cong &
\bigoplus_{i = 0}^{d} 
L_{n-i}^{\langle -\infty \rangle}(R)^{\binom{d}{i}}.
\end{eqnarray*}
\end{theorem}

\begin{example}[$\IZ^d \times G$] \label{exa:Zd_times_Z/p}
Let $G$ be a  group. By Theorem~\ref{the:Zd},
\begin{multline*}
K_n(R[G \times \IZ^d]) \cong K_n(RG[\IZ^d])
\\
\cong \bigoplus_{i = 0}^{d} K_{n-i}(RG)^{\binom{d}{i}}  \oplus
\bigoplus_{C \in \calmicyc(\IZ^d)} \; \bigoplus_{i = 0}^{d-1} 
\bigl(\NK_{n-i}(RG) \oplus \NK_{n-i}(RG)\bigr)^{\binom{d-1}{i}},
\end{multline*}
where $\calmicyc(\IZ^d)$ is the set of maximal infinite cyclic subgroups of
$\IZ^d$. Since
$$H_n\bigl(B(G \times \IZ^d);\bfK(R)\bigr) 
\cong 
\bigoplus_{i = 0}^{d} H_n\bigl(BG;\bfK(R)\bigr)^{\binom{d}{i}},$$
this implies
\begin{eqnarray*}
\Wh_n(G \times \IZ^d;R) 
& \cong &
\bigoplus_{i = 0}^{d} \Wh_{n-i}(G;R)^{\binom{d}{i}} 
\\
& & \quad \oplus \bigoplus_{C \in \calmicyc(\IZ^d)} \bigoplus_{i = 0}^{d-1} \;
\bigl(\NK_{n-i}(RG) \oplus \NK_{n-i}(RG)\bigr)^{\binom{d-1}{i}}.
\end{eqnarray*}
\end{example}

\begin{example}[Surface groups] \label{exa:surface_groups}
Let $\Gamma_g$ be the fundamental group of the orientable closed
surface of genus $g$, and let $\calm$ be a complete system of representatives of 
the conjugacy classes of maximal infinite cyclic subgroups of $G$.
If $g = 0$, then $\Gamma_g$ is trivial. If
$g = 1$, then $\Gamma_g$ is $\IZ^2$ and
Theorem~\ref{the:Zd} implies
\[
K_n(R\Gamma_1) 
\; \cong \;  
K_n(R) \oplus K_{n-1}(R)^2  \oplus K_{n-2}(R)
\oplus \bigoplus_{V \in \calm}  \bigl(\NK_n(R)^2  \oplus \NK_{n-1}(R)^2\bigr)\]
and
\[
L_n^{\langle - \infty \rangle}(R\Gamma_1) 
\; \cong \;
L^{\langle - \infty \rangle}_n(R) 
\oplus L_{n-1}^{\langle - \infty  \rangle}(R)^{2} 
\oplus L^{\langle - \infty \rangle}_{n-2}(R).\]
If $g \ge 2$, then $\Gamma_g$ is hyperbolic and torsion-free, so by
Theorem~\ref{the:torsionfree_hyperbolic_groups}
\[
K_n(R\Gamma_g) 
\; \cong \; 
H_{n}\bigl(B\Gamma_g;\bfK(R))\oplus  
\bigoplus_{V \in   \calm}  \NK_n(R)^2,\]
and
\[
L_n^{\langle - \infty \rangle}(R\Gamma_g) 
\; \cong \;
H_{n}\bigl(B\Gamma_g;\bfL^{\langle - \infty  \rangle}(R)\bigr).
\]
Since $\Gamma_g$ is stably a product of spheres,
\[
H_{n}\bigl(B\Gamma_g;\bfK(R))
\; \cong \; 
K_n(R) \oplus K_{n-1}(R)^{2g}\oplus K_{n-2}(R),\]
and
\[
H_{n}\bigl(B\Gamma_g;\bfL^{\langle - \infty  \rangle}(R)\bigr)
\; \cong \; 
L^{\langle - \infty \rangle}_n(R) \oplus L^{\langle - \infty \rangle}_{n-1}(R)^{2g}\oplus L^{\langle - \infty \rangle}_{n-2}(R).
\]

If $R = \IZ$, then for every $i \in \{1,0, -1, \ldots \} \amalg \{-\infty\}$, 
$L_n^{\langle i \rangle}(\IZ)$ is $\IZ$ if $n \equiv 0 \mod 4$,
$\IZ/2$ if $n \equiv 2 \mod  4$, and is trivial otherwise. Therefore,
$$L_n^{\langle i \rangle}(\IZ \Gamma_g) 
\cong
\begin{cases}\IZ \oplus \IZ/2 &  \text{if} \; n \equiv 0,2 \!\mod 4;
\\
\IZ^g &  \text{if} \; n \equiv 1 \! \mod 4;\\
(\IZ/2)^g &   \text{if} \;  n \equiv 3 \! \mod  4.
\end{cases}
$$
\end{example}

More generally, cocompact planar groups 
(sometimes called cocompact non-Euclidean crystallographic groups),
e.g., cocompact Fuchsian groups,
are treated in L\"uck-Stamm~\cite{Lueck-Stamm(2000)}
for $R = \IZ$. These computations can be carried over to arbitrary $R$.


\subsection{Virtually finitely generated abelian groups}
\label{subsec:K-theory_in_the_case_of_a_free_conjugation_action}


Consider the group extension
\begin{eqnarray}
& 1 \to A \to G \xrightarrow{q} Q \to 1,
\label{1_to_A_to_G_to_Q_to_1}
\end{eqnarray}
where $A$ is isomorphic to $\IZ^d$ for some $d \ge 0$ and $Q$ is a finite
group. The conjugation action of $G$ on the normal abelian subgroup $A$ induces
an action of $Q$ on $A$ via a group homomorphism which we denote by $\rho \colon
Q \to \aut(A)$.

Let $\calmicyc(A)$ be the set of maximal infinite cyclic subgroups of
$A$.  Since any automorphism of $A$ sends a maximal infinite cyclic
subgroup to a maximal infinite cyclic subgroup,
$\rho$ induces a $Q$-action on $\calmicyc(A)$.  Fix a
subset
\begin{eqnarray*}
I & \subseteq & \calmicyc(A)
\end{eqnarray*}
such that the intersection of every $Q$-orbit in $\calmicyc(A)$ with
$I$ consists of precisely one element.

For $C \in I$, let
\begin{eqnarray*}
Q_C & \subseteq & Q
\end{eqnarray*}
be the isotropy group of $C \in \calmicyc(A)$ under the $Q$-action. Let
\begin{eqnarray*}
I_1 & = & \bigl\{C \in I \mid Q_C = \{1\}\bigr\}, 
\\
I_2 & = & \bigl\{C \in I \mid Q_C = \IZ/2\bigr\}, 
\end{eqnarray*}
and let $J$ be a complete system of representatives of maximal non-trivial
finite subgroups of $G$.


\subsubsection{$K$-theory in the case of a free conjugation action}
\label{subsubsec:K-theory_in_the_case_of_a_free_conjugation_action}


\begin{theorem}
\label{the:K-theory_of_G_with_Aq_is_0}
Consider the group extension $ 1 \to A \to G \xrightarrow{q} Q \to 1$,
where $A$ is isomorphic to $\IZ^d$ for some $d \ge 0$ and $Q$ is a
finite group.  Suppose that the $Q$-action on $A$ is free away from $0 \in A$.

\begin{enumerate}

\item \label{the:K-theory_of_G_with_Aq_is_0:eimplicite}
For each $n \in \IZ$ there is an isomorphism induced by the various inclusions
\[
\left(\bigoplus_{F \in J} \Wh_n(F;R)\right) 
\oplus 
\bigl(\IZ \otimes_{\IZ Q} \Wh_n(A;R)\bigr) \xrightarrow{\cong} 
\Wh_n(G;R).
\]

\item \label{the:K-theory_of_G_with_Aq_is_0:explicite}

For each integer $n$, $\IZ \otimes_{\IZ Q} \Wh_n(A;R)$ is isomorphic to
$$
\left(\bigoplus_{C \in I_1} 
\bigoplus_{i=0}^{d-1} \bigl(\NK_{n-i}(R) \oplus \NK_{n-i}(R)\bigr)^{\binom{d-1}{i}}
\right) 
\\
\oplus \left(\bigoplus_{C \in I_2} 
\bigoplus_{i = 0} ^{d-1}  \NK_{n-i}(R)^{\binom{d-1}{i}}\right).$$

\end{enumerate}
\end{theorem}

If $Q = \{1\}$, then Theorem~\ref{the:K-theory_of_G_with_Aq_is_0} reduces to
Theorem~\ref{the:Zd} since $J = \emptyset$, $I_2 = \emptyset$ and $I_1$ is the
set of maximal infinite cyclic subgroups of $G = \IZ^d$.

\begin{theorem} \label{the:K-theory_for_regular_R_and_R_is_Z}
Under the assumptions of Theorem~\ref{the:K-theory_of_G_with_Aq_is_0}:
\begin{enumerate}
\item \label{the:K-theory_for_regular_R_and_R_is_Z:R_regular}
If $R$ is regular, then for each $n \in \IZ$
$$\bigoplus_{F \in J} \Wh_n(F;R) \cong  \Wh_n(G;R).$$
In particular,
$$\bigoplus_{F \in J} K_n(RF) \cong  K_n(RG) \quad \text{for} \; n \le -1$$
and
$$\bigoplus_{F \in J} \coker\bigl(K_0(R) \to K_0(RF)\bigr) \cong  
\coker\bigl(K_0(R) \to K_0(RG)\bigr).$$

\item \label{the:K-theory_for_regular_R_and_R_is_Z:R_Dedekind}
If $R$ is a Dedekind ring  of characteristic zero, then 
$$K_n(RG) = \{0\}  \quad \text{for} \; n \le -2.$$
\end{enumerate}
\end{theorem}


Applying Theorem~\ref{the:K-theory_for_regular_R_and_R_is_Z} to the special case
$R=\IZ$, we recover the fact that if $G$ satisfies the assumptions of
Theorem~\ref{the:K-theory_of_G_with_Aq_is_0}, then:
\begin{eqnarray*}
\bigoplus_{F \in J} \Wh(F) & \xrightarrow{\cong} & \Wh(G);
\\
\bigoplus_{F \in J} \widetilde{K}_0(\IZ F)  & \xrightarrow{\cong}  & \widetilde{K}_0(\IZ G);
\\
\bigoplus_{F \in J} K_{-1}(\IZ F)  & \xrightarrow{\cong}  & K_{-1}(\IZ G);
\\
K_n(\IZ G)  & \cong & \{0\} \quad \text{for} \; n \le -2.
\end{eqnarray*}


\subsubsection{$L$-theory in the case of a free conjugation action}
\label{subsubsec:L-theory_in_the_case_of_a_free_conjugation_action}

\begin{theorem}
\label{the:L-theory_of_G_with_Aq_is_0}
Consider the group extension $ 1 \to A \to G \xrightarrow{q} Q \to 1$, where
$A$ is isomorphic to $\IZ^d$ for some $d \ge 0$ and $Q$ is a finite group.
Suppose that the $Q$-action on $A$ is free away from $0 \in A$.  Assume that
there exists $n_0 \le -2$ such that $K_n(RV) = 0$ for all $n \le n_0$ and all
virtually cyclic subgroups $V\subseteq G$.  (By
Theorem~\ref{the:K-theory_for_regular_R_and_R_is_Z}~\ref{the:K-theory_for_regular_R_and_R_is_Z:R_Dedekind}, 
this condition is satisfied if $R$ is a Dedekind ring of characteristic zero).

Then there is an isomorphism
$$
\left(\bigoplus_{F \in J} \cals_n^{\per,\langle - \infty \rangle}(F;R)\right)
\oplus
\left(\bigoplus_{C \in I_2} \bigoplus_{H \in J_C} 
\UNil^{\langle - \infty \rangle}_n(D_{\infty};R)\right)
\xrightarrow{\cong} 
\cals_n^{\per,\langle -\infty \rangle}(G;R),
$$
where $J_C$ is a complete system of representatives of the conjugacy classes of
maximal finite subgroups of $W_GC=N_GC/C$. 

If $Q$ has odd order, then

$$
\bigoplus_{F \in J} \cals_n^{\per,\langle - \infty \rangle}(F;R)
\xrightarrow{\cong} 
\cals_n^{\per,\langle -\infty \rangle}(G;R).
$$
\end{theorem}

Dealing with groups extensions of the type described above when the conjugation
action of $Q$ on $A$ is not free is considerably harder than the free
case. However, we are able to say something when $Q$ is a cyclic group of prime
order and $R$ is regular.


\subsubsection{$K$-theory in the case $Q = \IZ/p$ for a prime $p$ and regular $R$}
\label{subsubsec:K-theory_in_the_case_Q_is_Z/p_for_a_prime_p}

\begin{theorem}
\label{the:K-theory_RG_Q_is_Z/p}
Consider the group extension $ 1 \to A \to G \xrightarrow{q} \IZ/p \to 1$,
where $A$ is isomorphic to $\IZ^d$ for some $d \ge 0$ and $p$ is a
prime number. Let $e$ be the natural number given by $A^{\IZ/p} \cong \IZ^e$, $J$ be
a complete system of representatives of the conjugacy classes of
non-trivial finite subgroups of $G$, and $\calmicyc(A^{\IZ/p})$ be the set of
maximal infinite cyclic subgroups of $A^{\IZ/p}$.

If $R$ is regular, then there is an isomorphism
\begin{multline*}
\Wh_n(G;R) \cong \left( \bigoplus_{H\in J} 
\bigoplus_{i=0}^e  \Wh_{n-i}(H;R)^{\binom{e}{i}}\right)
\\
\oplus
\left( \bigoplus_{H\in J}\; \bigoplus_{C \in  \calmicyc(A^{\IZ/p})} \;
\bigoplus_{i = 0}^{e-1} 
\bigl(\NK_{n-i}(R[\IZ/p]) \oplus \NK_{n-i}(R[\IZ/p])\bigr)^{\binom{e-1}{i}}\right).
\end{multline*}
\end{theorem}

\begin{remark}[Cardinality of $J$]\label{rem:cardinality_of_J}
Assume that the group $G$ appearing in
Theorem~\ref{the:K-theory_RG_Q_is_Z/p}
is not torsion-free, or, equivalently, that $J$ is non-empty. 
Consider $A$ as a $\IZ[\IZ/p]$-module via
the conjugation action $\rho \colon \IZ/p \to \aut(A)$.
Then there is a bijection
\begin{eqnarray*}
H^1(\IZ/p;A) & \xrightarrow{\cong} &  J,
\end{eqnarray*}
defined as follows. Fix an element $t \in G$ of order $p$.
Every $\overline{x} \in H^1(\IZ/p;A)$ is represented by an
element $x$ in the kernel of $\sum_{i=0}^{p-1} \rho^i \colon A \to
A$. Then $xt$ has order $p$. Send $\overline{x}$ to the unique element of $J$ 
that is conjugate to $\langle xt \rangle$ in $G$.
\end{remark}

If $p$ is prime and $R$ is regular, 
then Example~\ref{exa:Zd_times_Z/p} in the case $G = \IZ/p$  
is consistent with Theorem~\ref{the:K-theory_RG_Q_is_Z/p} and
Remark~\ref{rem:cardinality_of_J}. Namely,  $\IZ/p$ acts trivially on
$\IZ^d$, and so $A = A^{\IZ/p}$, $J=\{\{0\} \times \IZ/p\}$, and $d = e$.


\subsubsection{$L$-theory in the case $Q = \IZ/p$ for an odd prime $p$.}
\label{subsubsec:L-theory_in_the_case_Q_is_Z/p_for_an_odd_prime_p}

\begin{theorem}
\label{the:L-theory_RG_Q_is_Z/p}
Consider the group extension $ 1 \to A \to G \xrightarrow{q} \IZ/p \to 1$, 
where $A$ is isomorphic to $\IZ^d$ for some $d \ge 0$ and
$p$ is an odd prime number.  Let $e$ be the natural number given by
$A^{\IZ/p} \cong \IZ^e$ and $J$ be a complete system of
representatives of the conjugacy classes of non-trivial finite
subgroups of $G$.

Then there is an isomorphism
$$\bigoplus_{H\in J} 
\bigoplus_{i=0}^e  \cals^{\per,\langle - \infty \rangle}_{n-i}(H;R)^{\binom{e}{i}}
\xrightarrow{\cong} 
\cals^{\per,\langle - \infty \rangle}_{n}(G;R).
$$
\end{theorem}

Two-dimensional crystallographic groups are treated in
Pearson~\cite{Pearson(1998)} and L\"uck-Stamm~\cite{Lueck-Stamm(2000)}
for $R = \IZ$. Some of the computations there can be carried over to arbitrary $R$.
The Whitehead groups of three-dimensional crystallographic groups
are computed in Alves-Ontaneda~\cite{Alves-Ontaneda(2006)}.

In the case $R = \IZ$, we get a computation for all decorations.

\begin{theorem}
\label{the:L-theory_ZG_Q_is_Z/p}
Consider the group extension $ 1 \to A \to G \xrightarrow{q} \IZ/p \to 1$, 
where $A$ is isomorphic to $\IZ^d$ for some $d \ge 0$ and
$p$ is an odd prime number.  Let $e$ be the natural number given by
$A^{\IZ/p} \cong \IZ^e$, $J$ be a complete system of
representatives of the conjugacy classes of non-trivial finite
subgroups of $G$, and $\epsilon$ be any of the decorations $s$, $h$, $p$ or
$\{\langle j \rangle \mid j = 1,0,-1, \ldots\} \amalg \{-\infty\}$.
Define
\begin{eqnarray*}
\cals_n^{\per,\epsilon}(G;\IZ) 
& := &
H_n^G\bigl(EG \to \pt;\bfL_{\IZ}^{\epsilon}\bigr).
\end{eqnarray*}

Then there is an isomorphism
$$\bigoplus_{H\in J} 
\bigoplus_{i=0}^e  \cals^{\per,\epsilon}_{n-i}(H;\IZ)^{\binom{e}{i}}
\xrightarrow{\cong} 
\cals^{\per,\epsilon}_{n}(G;\IZ).
$$
\end{theorem}

\begin{remark} In general, the structure sets $\cals^{\per,\epsilon}_{n-i}(H;\IZ)$ depend 
on the decoration $\epsilon$. We mention without proof  that  for an odd prime $p$
$$\cals_n^{\per,s}(\IZ/p;\IZ) 
\cong \widetilde{L}^s_n(\IZ[\IZ/p])[1/p]
\cong \begin{cases}
\IZ[1/p]^{(p-1)/2} & n \; \text{even};
\\
\{0\} & n \; \text{odd}.
\end{cases}
$$
\end{remark}

\begin{remark} The computation of the structure set
$\cals^{\per,\epsilon}_{n}(G;\IZ)$ when the conjugation
action of $\IZ/p$ on $\IZ^d$ is free plays a role in a forthcoming
paper by Davis and L\"uck, in which this case is further analyzed to compute the geometric structure sets of certain
manifolds that occur as total spaces of a bundle over lens
spaces with $d$-dimensional tori as fibers.

\end{remark}

\begin{remark}[Topological $K$-theory of reduced group $C^*$-algebras]
\label{rem:Topological_K-theory_of_reduced_group_Cast-algebras}
All of the above computations also apply to the topological $K$-theory of the
reduced group $C^*$-algebra, since the Baum-Connes Conjecture is true for
these groups.  In the Baum-Connes setting the situation simplifies
considerably because one works with $\eub{G}$ instead of $\edub{G}$ and hence
there are no $\Nil$-phenomena. For instance, if $G$ is a hyperbolic group, then there is an
isomorphism
$$K_n^G(\eub{G}) \xrightarrow{\cong} K_n(C^*_r(G))$$
from the equivariant topological $K$-theory of $\eub{G}$ to the topological
$K$-theory of the reduced group $C^*$-algebra $C^*_r(G)$. In the case of a
torsion-free hyperbolic group, this reduces to an isomorphism
$$K_n(BG) \xrightarrow{\cong} K_n(C^*_r(G)).$$
The case $G = \IZ^d \rtimes_{\rho} \IZ/p$ for a free conjugation action has been
carried out for both complex and real topological $K$-theory in detail
in~\cite[Theorem~0.3 and Theorem~0.6]{Davis-Lueck(2010)}.
\end{remark}


\typeout{--------------------------------------   Section 2 :  Background  ----------------------------------}

\section{Background}\label{sec:Background}

In this section we give some background about the Farrell-Jones Conjecture and related topics.


\subsection{Classifying Spaces for Families}\label{subsec:Classifying_Spaces_for_Families}

Let $G$ be a group. A \emph{family of subgroups} of $G$ is a
collection of subgroups that is closed under conjugation and taking
subgroups. Examples of such families are
\begin{align*}
\caltr    & = \{\text{trivial subgroup}\}; \\
\calfin   & = \{\text{finite subgroups}\}; \\
\calvcyc  &= \{\text{virtually cyclic subgroups}\};\\
\calall &= \{\text{all subgroups}\}.
\end{align*}

Let $\calf$ be a family of subgroups of $G$. A model for the
\emph{universal space $\EGF{G}{\calf}$ for $\calf$} is a
$G$-$CW$-complex $X$ with isotropy groups in $\calf$ such that for any
$G$-$CW$-complex $Y$ with isotropy groups in $\calf$ there exists a
$G$-map $Y \to X$ that is unique up to $G$-homotopy.  In other words,
$X$ is a terminal object in the $G$-homotopy category of
$G$-$CW$-complexes whose isotropy groups belong to $\calf$. In
particular, any two models for $\EGF{G}{\calf}$ are $G$-homotopy
equivalent, and for two families $\calf_0 \subseteq \calf_1$, there is
precisely one $G$-map $\EGF{G}{\calf_0} \to \EGF{G}{\calf_1}$ up to
$G$-homotopy.

For every group $G$ and every family of subgroups $\calf$ there
exists a model for $\EGF{G}{\calf}$. A $G$-$CW$-complex $X$ is a model
for $\EGF{G}{\calf}$ if and only if the $H$-fixed point set $X^H$ is
contractible for every $H$ in $\calf$ and is empty otherwise. For
example, a model for $\EGF{G}{\calall}$ is $G/G=\pt$, and a model for
$\EGF{G}{\caltr}$ is the same as a model for $EG$, the total space of
the universal $G$-principal bundle $EG \to BG$. The \emph{universal
$G$-$CW$-complex for proper $G$-actions}, $ \EGF{G}{\calfin}$, will
be denoted by $\eub{G}$, and the universal space $\EGF{G}{\calvcyc}$
for $\calvcyc$ will be denoted by $\edub{G}$. For more information on
classifying spaces the reader is referred to the survey article by
L\"uck~\cite{Lueck(2005s)}.


\subsection{Review of the Farrell-Jones
Conjecture}\label{subsec:Review_of_the_Farrell-Jones_Conjecture}
Let $\calh^?_*$ be an equivariant homology theory in the sense of
L\"uck~\cite[Section~1]{Lueck(2002b)}. 
Then, for every group
$G$ and every $G$-$CW$-pair $(X,A)$ there is a $\IZ$-graded
abelian group $\calh^G_*(X,A)$, and subsequently a $G$-homology theory $\calh^G_*$.
For every group homomorphism $\alpha \colon H \to G$, every
$H$-$CW$-pair $(X,A)$ and every $n \in \IZ$, there is a natural homomorphism
$\ind_\alpha:\calh^H_*(X,A)\to \calh^G_*(G \times_{\alpha}(X,A))$, known as the
induction homomorphism. If the kernel of $\alpha$ operates relative freely on $(X,A)$, then $\ind_\alpha$ is an
isomorphism.

Our main examples are the equivariant homology theories
$H_*^?(-;\bfK_R)$ and $H_*^?(-;\bfL_R^{\langle - \infty
\rangle})$ appearing in the $K$-theoretic and $L$-theoretic
Farrell-Jones Conjectures, where $R$ is an associative ring with unit
(and involution) (see
L\"uck-Reich~\cite[Section~6]{Lueck-Reich(2005)}).  The basic property
of these two equivariant homology theories is that
$$
\begin{array}{lclcl}
H_n^G\bigl(G/H;\bfK_R\bigr) 
& = & 
H_n^H\bigl(\pt;\bfK_R\bigr) 
& = & 
K_n(RH);
\\
H_n^G\bigl(G/H;\bfL_R^{\langle - \infty \rangle}\bigr) 
& = & 
H_n^H\bigl(\pt;\bfL_R^{\langle - \infty \rangle}\bigr) 
& = & 
L_n^{\langle - \infty \rangle}(RH),
\end{array}
$$
for every subgroup $H \subseteq G$.

The \emph{Farrell-Jones Conjecture for a group $G$ and a ring $R$}
states that the assembly maps
\begin{eqnarray*}
H_n^G\bigl(\edub{G};\bfK_R\bigr) 
& \xrightarrow{\cong} & 
H_n^G\bigl(G/G;\bfK_R\bigr) = K_n(RG);
\\
H_n^G\bigl(\edub{G};\bfL_R^{\langle - \infty \rangle}\bigr) 
& \xrightarrow{\cong} & 
H_n^G\bigl(G/G;\bfL_R^{\langle - \infty \rangle}\bigr) = L_n^{\langle - \infty \rangle}(RG),
\end{eqnarray*}
induced by the projection $\edub{G} \to G/G$, are bijective for every
$n \in \IZ$~\cite{Farrell-Jones(1993a)}. The Farrell-Jones Conjecture
has been studied extensively because of its geometric significance. It
implies, in dimensions $\ge 5$, both the Novikov Conjecture about the homotopy invariance of higher signatures
and the Borel Conjecture about the rigidity of
manifolds with fundamental group $G$. It also implies other well-known
conjectures, such as the ones due to Bass and Kadison. For a survey and
applications of the Farrell-Jones Conjecture, see, for example,
L\"uck-Reich~\cite[Section~6]{Lueck-Reich(2005)} and
Bartels-L\"uck-Reich~\cite{Bartels-Lueck-Reich(2008appl)}.

\begin{theorem} [Farrell-Jones Conjecture for hyperbolic groups and
virtually $\IZ^d$-groups]
\label{the:FJC_hyperbolic_virtually_Zd}
The Farrell-Jones Conjecture is true if $G$ is a hyperbolic group or
a virtually finitely generated abelian group, and $R$ is any
ring. 
\end{theorem}

The $K$-theoretic Farrell-Jones Conjecture with coefficients in any
additive $G$-category for a hyperbolic group $G$ was proved by
Bartels-Reich-L\"uck~\cite{Bartels-Lueck-Reich(2008hyper)}. The
version of the Farrell-Jones Conjecture with coefficients in an
additive category encompasses the version with rings as
coefficients. The $L$-theoretic Farrell-Jones Conjecture with
coefficients in any additive $G$-category for hyperbolic groups and
$\CAT(0)$-groups was established by Bartels and
L\"uck~\cite{Bartels-Lueck(2012annals)}. In that paper it is also
shown that the $K$-theoretic assembly map is $1$-connected for
$\CAT(0)$-groups. Note that a virtually finitely generated abelian
group is $\CAT(0)$. Quinn~\cite[Theorem~1.2.2]{Quinn(2005)}) proved
that the $K$-theoretic assembly map for virtually finitely generated
abelian groups is bijective for every integer $n$ if $R$ is a
commutative ring. However, the proof carries over to the
non-commutative setting.

\begin{remark}[The Interplay of $K$- and $L$-Theory]
\label{rem:Rothenberg_sequence} $L$-theory $L_*^{\langle i
  \rangle}(RG)$ can have various decorations for $i \in \{2,1,0,-1,-2,
\ldots\} \amalg \{-\infty\}$. One also finds $L_*^{\epsilon}(RG)$
for $\epsilon = p,h,s$ in the literature.  The decoration $\langle 1
\rangle $ coincides with the decoration $h$, $\langle 0 \rangle $
with the decoration $p$, and $\langle 2 \rangle $ is related to the
decoration $s$.  For $j \leq 1$ there are forgetful maps
$L_n^{\langle j +1 \rangle}(R) \to L_n^{\langle j \rangle}(R)$ that
fit into the so-called \emph{Rothenberg sequence} (see
Ranicki~\cite[Proposition~1.10.1 on page~104]{Ranicki(1981)},
\cite[17.2]{Ranicki(1992a)})
\begin{multline}
  \cdots \to L_n^{\langle j+1 \rangle}(R ) \to L_n^{\langle j
    \rangle}(R ) \to \widehat{H}^n(\IZ/ 2; \widetilde{K}_{j}(R)) \\
  \to L_{n-1}^{\langle j+1 \rangle}(R ) \to L_{n-1}^{\langle j
    \rangle}(R) \to \cdots.
  \label{Rothenberg_sequence}
\end{multline}
$\widehat{H}^n(\IZ/ 2;\widetilde{K}_{j}(R))$ denotes the
Tate-cohomology of the group $\IZ/ 2$ with coefficients in the
$\IZ[\IZ/ 2]$-module $\widetilde{K}_{j}(R)$. The involution on
$\widetilde{K}_{j}(R)$ comes from the involution on $R$.  There is a
similar sequence relating $L^s_n (RG)$ and $L^h_n (RG)$, where the
third term is the $\IZ/2$-Tate-cohomology of $\Wh^R_1 (G)$.

For geometric applications the most important case is $R=\IZ$ with
the decoration $s$. If $\widetilde{K}_0(\IZ G)$, $\Wh(G)$ and
$K_n(\IZ G)$ for $n \le -1$ all vanish, then the Rothenberg sequence
for $n \in \IZ$ implies the bijectivity of the natural maps

$$L_n^{s}(\IZ G) \xrightarrow{\cong} 
L_n^{h}(\IZ G) \xrightarrow{\cong} L_n^{p}(\IZ G) \xrightarrow{\cong}
L_n^{\langle - \infty \rangle}(\IZ G).$$

In the formulation of the Farrell-Jones Conjecture (see
Section~\ref{sec:Background}), one must use the decoration $\langle -
\infty \rangle$ since the conjecture is false otherwise (see
Farrell-Jones-L\"uck~\cite{Farrell-Jones-Lueck(2002)}).
\end{remark}


\subsection{Equivariant homology and relative assembly maps}
\label{subsec:Equivariant_homology_and_relative_assembly_maps}

Consider an equivariant homology theory $\calh^?_*$ in the sense
of~\cite[Section~1]{Lueck(2002b)}.  Given a $G$-map $f \colon X \to Y$
of $G$-$CW$-complexes and $n \in \IZ$, define
\begin{eqnarray}
\calh_n^G(f) &:= & \calh_n^G\bigl(\cyl(f_0),X\bigr)
\label{calh_nG(f)}
\end{eqnarray}
for any cellular $G$-map $f_0 \colon X \to Y$ that is $G$-homotopic to
$f$.  Here, $\cyl(f_0)$ is the $G$-$CW$-complex given by the mapping
cylinder of $f$. It contains $X$ as a $G$-$CW$-subcomplex.  Such an
$f_0$ exists by the Equivariant Cellular Approximation Theorem
(see~\cite[Theorem~II.2.1 on page~104]{Dieck(1987)}).  The definition
is independent of the choice of $f_0$, since two cellular
$G$-homotopic $G$-maps $f_0,f_1 \colon X \to Y$ are cellularly
$G$-homotopic by the Equivariant Cellular Approximation Theorem
(see~\cite[Theorem~II.2.1 on page~104]{Dieck(1987)}) and a cellular
$G$-homotopy between $f_0$ and $f_1$ yields a $G$-homotopy equivalence
$u \colon \cyl(f_0) \to \cyl(f_1)$ that is the identity on $X$. This
implies that $\calh_n^G(f)$ depends only on the $G$-homotopy class of
$f$. From the axioms of an equivariant homology theory, $\calh_n^G(f)$
fits into a long exact sequence
\begin{eqnarray}
\cdots \to \calh_{n+1}^G(f) \to \calh_n^G(X) \to \calh_n^G(Y) \to
\calh_n^G(f) \to \calh_{n-1}^G(X) \to \cdots. 
\label{long_exact}
\end{eqnarray}

The following fact is proved in Bartels~\cite{Bartels(2003b)}.

\begin{lemma}
\label{lem:splitting_eub_deub}
For every group $G$, every ring $R$, and
every $n \in \IZ$:
\begin{enumerate} 
\item \label{lem:splitting_eub_deub:K}
the relative assembly map
$$H_n^G\bigl(\eub{G};\bfK_R\bigr)  \to  H_n^G\bigl(\edub{G};\bfK_R\bigr)$$
is split-injective;

\item \label{lem:splitting_eub_deub:L}
the relative assembly map
$$H_n^G\bigl(\eub{G};\bfL_R^{\langle - \infty \rangle}\bigr)
\to
H_n^G\bigl(\edub{G};\bfL_R^{\langle - \infty \rangle}\bigr)$$
is split-injective, provided there is an $n_0 \le -2$ such that
$K_n(RV) = 0$ for every $n \le n_0$ and every virtually cyclic subgroup
$V \subseteq G$. 

\end{enumerate}
\end{lemma}

\begin{remark} \label{lem:splitting_eub_deub:con}
Notice that the condition appearing in assertion~\ref{lem:splitting_eub_deub:L}
of the above lemma is automatically satisfied if one of the following stronger conditions holds:
\begin{enumerate}
\item $R = \IZ$;
\item $R$ is regular with $\IQ \subseteq R$.
\end{enumerate}
If $R = \IZ$, this follows from~\cite[Theorem~2.1(a)]{Farrell-Jones(1995)}.
If $R$ is regular and $\IQ \subseteq R$,
then for every virtually cyclic subgroup $V \subseteq G$ the ring
$RV$ is regular and hence $K_n(RV) = 0$ for all $n \le -1$.
\end{remark}

Lemma~\ref{lem:splitting_eub_deub} tells us that the source of the
assembly map appearing in the Farrell-Jones Conjecture can be computed
in two steps, the computation of $H_n^G\bigl(\eub{G};\bfK_R\bigr)$ and
the computation of the remaining term
$H_n^G\bigl(\eub{G}\to\edub{G};\bfK_R\bigr)$ defined
in~\eqref{calh_nG(f)}. Furthermore, Lemma~\ref{lem:splitting_eub_deub}
implies
\begin{eqnarray}
H_n^G\bigl(EG\to \edub{G};\bfK_R\bigr) \cong 
H_n^G\bigl(EG\to \eub{G};\bfK_R\bigr) \oplus H_n^G\bigl(\eub{G}\to \edub{G};\bfK_R\bigr);
\label{relative_K_splitting}
\end{eqnarray}
and
\begin{multline}
H_n^G\bigl(EG\to \edub{G};\bfL_R^{\langle -\infty \rangle}\bigr) \cong \\
H_n^G\bigl(EG\to \eub{G};\bfL_R^{\langle -\infty \rangle}\bigr) \oplus H_n^G\bigl(\eub{G}\to \edub{G};\bfL_R^{\langle -\infty \rangle}\bigr),
\label{relative_L_splitting}
\end{multline}
provided that there exists an $n_0 \le -2$ such that $K_n(RV) = 0$ for
every $n \le n_0$ and every virtually cyclic subgroups $V \subseteq
G$. This is useful for calculating the Whitehead groups for a given
group $G$ and ring $R$ which satisfy the Farrell-Jones Conjecture.

\begin{remark}\label{K-theory_induction}
The induction structure of $H_*^?(-;\bfK_R)$ can be used to define a $\IZ/2$-action on $H_n^{\IZ}\bigl(E\IZ \to \pt;\bfK_R\bigr)$ that is compatible with the $\IZ/2$-action on $\NK_n(R) \oplus \NK_n(R)$ given by flipping the two factors. The action is defined by the composition of isomorphisms
\begin{multline*} 
\tau:H_n^{\IZ}\bigl(E\IZ \to \pt;\bfK_R\bigr) \xrightarrow{\ind_{-\id_{\IZ}}} H_n^{\IZ}\bigl(\ind_{-\id_{\IZ}}E\IZ \to \pt;\bfK_R\bigr)
\\ \xrightarrow{} H_n^{\IZ}\bigl(E\IZ \to \pt;\bfK_R\bigr), 
\end{multline*}
where the second map is induced by the unique (up to equivariant homotopy) equivariant map $l:\ind_{-\id_{\IZ}}E\IZ \to E\IZ$; $l$ is an equivariant homotopy equivalence since $\ind_{-\id_{\IZ}}E\IZ$ is a model for $E\IZ$. 

To see that this corresponds to the flip action on $\NK_n(R) \oplus \NK_n(R)$, consider the following diagram coming from the long exact sequence~\eqref{long_exact}. 
\[
\xymatrix{
	H_n^{\IZ}\bigl(E\IZ;\bfK_R\bigr)
		\ar[d]^{\ind_{-\id_{\IZ}}} \ar[r] &
	H_n^{\IZ}\bigl(\pt;\bfK_R\bigr)
		\ar[d]^{\ind_{-\id_{\IZ}}} \ar[r] &
	H_n^{\IZ}\bigl(E\IZ \to \pt;\bfK_R\bigr)
		\ar[d]^{\ind_{-\id_{\IZ}}}
		\\
	H_n^{\IZ}\bigl(\ind_{-\id_{\IZ}}E\IZ;\bfK_R\bigr)
		\ar[d]^{l_*} \ar[r] &
	H_n^{\IZ}\bigl(\pt;\bfK_R\bigr)
		\ar[d]^= \ar[r] &
	H_n^{\IZ}\bigl(\ind_{-\id_{\IZ}}E\IZ \to \pt;\bfK_R\bigr)
		\ar[d]^{l_*}
		\\
	H_n^{\IZ}\bigl(E\IZ;\bfK_R\bigr)
		\ar[r] &
	H_n^{\IZ}\bigl(\pt;\bfK_R\bigr)
		\ar[r] &
	H_n^{\IZ}\bigl(E\IZ \to \pt;\bfK_R\bigr)
}
\]
Recall that $H_n^{\IZ}\bigl(E\IZ;\bfK_R\bigr) \to H_n^{\IZ}\bigl(\pt;\bfK_R\bigr)\cong K_n(R\IZ)\cong K_n(R[t,t^{-1}])$ is a split injection, and so the Bass-Heller-Swan decomposition~\eqref{Bass-Heller-Swan} establishes the identification $H_n^{\IZ}\bigl(E\IZ \to \pt;\bfK_R\bigr)\cong \NK_n(R) \oplus \NK_n(R)$~\eqref{H_nZ(EZ_to_pt;bfK_R)_cong_NK_n(R)_oplus_NK_n(R)}. Also recall that the two copies of $\NK_n(R)$ appearing in the decomposition of $K_n(R[t,t^{-1}])$ come from the embeddings $R[t] \hookrightarrow R[t,t^{-1}]$ and $R[t^{-1}] \hookrightarrow R[t,t^{-1}]$. By the definition of the induction structure for $H_*^?(-;\bfK_R)$, 
$$\ind_{-\id_{\IZ}}:H_n^{\IZ}\bigl(\pt;\bfK_R\bigr)\to H_n^{\IZ}\bigl(\pt;\bfK_R\bigr)$$
corresponds to the homomorphism $K_n(R[t,t^{-1}])\to K_n(R[t,t^{-1}])$ induced by interchanging $t$ and $t^{-1}$ (see, for example,~\cite[Section~6]{Lueck-Reich(2005)}), which swaps the two copies of $\NK_n(R)$ in the decomposition of $K_n(R[t,t^{-1}])$. Therefore the above diagram implies that 
$\tau$ coincides with the flip action on $\NK_n(R) \oplus \NK_n(R)$. A proof of this fact can also be found in~\cite[Lemma 3.22]{Davis-Quinn-Reich(2011)}.
\end{remark}

\begin{remark}[Relative assembly and Nil-terms]
\label{rem:Term_calh_nN_GVeub(N_GV)_to_EW_GV;bkk_R)}
Let $V$ be an infinite virtually cyclic group. If $V$ is of type I,
then $V$ can be written as a semi-direct product $F \rtimes \IZ$,
and $H^V_*\bigl(\eub{V} \to \pt;\bfK_R\bigr)$ can be identified with
the non-connective version of Waldhausen's Nil-term associated to
this semi-direct product (see~\cite[Sections~9
and~10]{Bartels-Lueck(2006)}). If $V$ is of type {II}, then it can
be written as an amalgamated product $V_1 \ast_{V_0} V_2$ of finite
groups, where $V_0$ has index two in both $V_1$ and $V_2$. In this
case, $H^V_*\bigl(\eub{V} \to \pt;\bfK_R\bigr)$ can be identified
with the non-connective version of Waldhausen's Nil-term associated
to this amalgamated product~\cite{Bartels-Lueck(2006)}.  The
identifications come from the Five-Lemma and the fact that both
groups fit into the same long exact sequence associated to the
semi-direct product, or amalgamated product, respectively. This is
analogous to the $L$-theory case which is explained below.  If $R$
is regular and $\IQ \subseteq R$, e.g., $R$ is a field of
characteristic zero, then $H^V_n\bigl(\eub{V}
\to \pt;\bfK_R\bigr) = 0$ for every virtually cyclic group
$V$, and hence, for any group $G$ the map
$$H_n^G\bigl(\eub{G};\bfK_R\bigr) \xrightarrow{\cong} H_n^G\bigl(\edub{G};\bfK_R\bigr)$$
is bijective (see~\cite[Proposition~2.6]{Lueck-Reich(2005)}).
\end{remark}


Let $\UNil_n^h(R;R,R)$ denote the Cappell
UNil-groups~\cite{Cappell(1974c)} associated to the amalgamated
product $D_{\infty} = \IZ/2 \ast \IZ/2$.  It is a direct summand in
$L_n^{h}(R[D_{\infty}])$ and there is a Mayer-Vietoris sequence
\begin{multline*}
\cdots \to L_n^{h}(R) \to L_n^{h}(R[\IZ/2]) \oplus L_n^{h}(R[\IZ/2])
\to L_n^{h}(R[D_{\infty}])/\UNil_n(R;R,R)
\\
\to L_n^{h}(R) \to L_n^{h}(R[\IZ/2]) \oplus L_n^{\langle -\infty
  \rangle}(R[\IZ/2]) \to \cdots.
\end{multline*}
There is also a Mayer-Vietoris sequence that maps to the one above
which comes from the model for $\eub{D_{\infty}}$ given by the obvious
$D_{\infty}$-action on $\IR$:
\begin{multline*}
\cdots \to L_n^{h}(R) \to L_n^{h}(R[\IZ/2]) \oplus L_n^{h}(R[\IZ/2])
\to H_n^{D_{\infty}}\bigl(\eub{D_{\infty}};\bfL_R^{h}\bigr)
\\
\to L_n^{h}(R) \to L_n^{h}(R[\IZ/2]) \oplus L_n^h(R[\IZ/2]) \to
\cdots.
\end{multline*}
This implies
\begin{eqnarray*}
\UNil_n^h(R;R,R) & \cong & 
H_*^{D_{\infty}}\bigl(\eub{D_{\infty}}  
\to  \pt;\bfL_R^{h}\bigr).
\end{eqnarray*}

\begin{lemma} \label{lem:UNil(R;R,R)_and_UNil-infty(D_infty;R)} Assume
that $\widetilde{K}_i(R) \cong \widetilde{K}_i(R[\IZ/2]) \cong
\widetilde{K}_i(R[D_{\infty}]) = 0$ for all $i \le 0$.  (This
condition is satisfied, for example, if $R = \IZ$.)  Then:
$$\UNil_n^h(R;R,R) \cong \UNil^{\langle - \infty \rangle}_n(D_{\infty};R).$$
\end{lemma}
\begin{proof}
Using the Rothenberg sequences~\eqref{Rothenberg_sequence}, one
obtains natural isomorphisms
\begin{eqnarray*}
  L_n^{\langle -\infty \rangle}(R) & \cong & L_n^{h}(R);
  \\
  L_n^{\langle -\infty \rangle}(R[\IZ/2]) & \cong & L_n^{h}(R[\IZ/2]);
  \\
  L_n^{\langle -\infty \rangle}(R[D_{\infty}]) & \cong & L_n^{h}(R[D_{\infty}]).
\end{eqnarray*}
Thus, a comparison argument involving the Atiyah-Hirzebruch spectral
sequences shows that the obvious map
$$H_n^{D_{\infty}}\bigl(\eub{D_{\infty}}   \to  \pt;\bfL_R^{h}\bigr)
\xrightarrow{\cong} H_n^{D_{\infty}}\bigl(\eub{D_{\infty}} \to
\pt;\bfL_R^{\langle -\infty \rangle}\bigr)$$ is bijective for all $n
\in \IZ$.
\end{proof}

Cappell's $\UNil$-terms~\cite{Cappell(1974c)} have been further
investigated
in~\cite{Banagl-Ranicki(2006)},\cite{Connolly-Ranicki(2005)}
and~\cite{Connolly-Davis(2004)}. The Waldhausen Nil-terms have been
analyzed in~\cite{Davis-Khan-Ranicki(2011)}
and~\cite{Grunewald(2008Nil)}.  If $2$ is inverted, the situation in
$L$-theory simplifies. Namely, for every $n \in \IZ$ and
every virtually cyclic group $V$, $H^V_n\bigl(\eub{V} \to
\pt;\bfL_R^{\langle \infty \rangle}\bigr)[1/2] = 0$. Therefore the map
$$H_n^G\bigl(\eub{G};\bfL_R^{\langle -\infty \rangle}\bigr)[1/2] 
\xrightarrow{\cong} H_n^G\bigl(\edub{G};\bfL_R^{\langle -\infty
\rangle}\bigr)[1/2]$$ is an isomorphism for any group $G$, and the
decorations do not play a role
(see~\cite[Proposition~2.10]{Lueck-Reich(2005)}).

\begin{remark}[Role of type I and II]
\label{rem:role_of_type_I_and_II}
Let $\calvcyc_I$ be the family of subgroups that are either finite
or infinite virtually cyclic of type I, i.e., groups admitting an
epimorphism onto $\IZ$ with finite kernel.  Then the following maps
are bijections
\begin{eqnarray*}
  H_n^G(\EGF{G}{\calvcyc_I};\bfK_R)
  &\xrightarrow{\cong} &
  H_n^G(\edub{G};\bfK_R);
  \\
  H_n^G\bigl(\eub{G};\bfL_R^{\langle - \infty \rangle}\bigr) 
  &\xrightarrow{\cong} &
  H_n^G\bigl(\EGF{G}{\calvcyc_I};\bfL_R^{\langle - \infty \rangle}\bigr).
\end{eqnarray*}
For the $K$-theory case, see, for instance,~\cite{Davis-Khan-Ranicki(2011),Davis-Quinn-Reich(2011)}.
The $L$-theory case is proven in
L\"uck~\cite[Lemma~4.2]{Lueck(2005heis)}).  In particular, for a
torsion-free group $G$, the map
\begin{eqnarray*}
  H_n^G\bigl(EG;\bfL_R^{\langle - \infty \rangle}\bigr) 
  &\xrightarrow{\cong} &
  H_n^G\bigl(\edub{G};\bfL_R^{\langle - \infty \rangle}\bigr)
\end{eqnarray*}
is a bijection.
\end{remark}


\typeout{------------------------------------------- Section 3: Proofs  -----------------------------------}

\section{Proofs of Results}\label{sec:proofs}
We now prove the results stated in Section~\ref{sec:statements}.


\subsection{Hyperbolic groups}
\label{subsec:proof_hyperbolic_groups}

\begin{proof}[Proof of Theorem~\ref{the:hyperbolic_groups}]
By~\cite[Corollary~2.11, Theorem~3.1 and Example~3.5]{Lueck-Weiermann(2012)} 
there is a $G$-pushout
$$
 \xymatrix{
 \coprod_{V\in\calm}G\times_{V}\eub{V}
 \ar[d]^{\coprod_{V\in\calm} p} \ar[r]^-i & \eub{G} \ar[d] \\
 \coprod_{V\in\calm}G/V \ar[r] & \edub{G}
 }
$$
where $i$ is an inclusion of $G$-$CW$-complexes, $p$ is the obvious
projection and $\calm$ is a complete system of representatives of the
conjugacy classes of maximal infinite virtually cyclic subgroups
of $G$.  Now Theorem~\ref{the:hyperbolic_groups} follows from
Theorem~\ref{the:FJC_hyperbolic_virtually_Zd} and
Lemma~\ref{lem:splitting_eub_deub}.~
\end{proof}

\begin{proof}[Proof of Theorem~\ref{the:torsionfree_hyperbolic_groups}]~%
\ref{the:torsionfree_hyperbolic_groups:K} 
This follows from~\eqref{H_nZ(EZ_to_pt;bfK_R)_cong_NK_n(R)_oplus_NK_n(R)} and
Theorem~\ref{the:hyperbolic_groups}~\ref{the:hyperbolic_groups:K}.
\\[1mm]~%
\ref{the:torsionfree_hyperbolic_groups:L}
Since any virtually cyclic subgroup of $G$ is trivial or infinite cyclic,
the claim follows from Theorem~\ref{the:FJC_hyperbolic_virtually_Zd} and
Remark~\ref{rem:role_of_type_I_and_II}.
\end{proof}


\subsection{$K$- and $L$-theory of $\IZ^d$}
\label{subsec:proof_Zd}

As a warm-up for virtually free abelian groups, we compute the $K$ and $L$-theory of $R[\IZ^d]$.

\begin{proof}[Proof of Theorem~\ref{the:Zd}]
Using the induction structure of the equivariant homology theory
$H_*^?(-;\bfK_R)$ and the fact that $B\IZ^d$ is the $d$-dimensional torus, it
follows that
\begin{eqnarray}
H^{\IZ^d}_n\bigl(E\IZ^d;\bfK_R\bigr) 
& \cong &
H^{\{1\}}_n\bigl(B\IZ^d;\bfK_R\bigr)
\nonumber
\\
& \cong & 
\bigoplus_{i = 0}^d H^{\{1\}}_{n-i}\bigl(\pt;\bfK_R\bigr)^{\binom{d}{i}}
\nonumber
\\
& = & 
\bigoplus_{i = 0}^d K_{n-i}(R)^{\binom{d}{i}}.
\label{HZd(EZd;bfk_R)}
\end{eqnarray}
Similarly one shows that
\begin{eqnarray}
H_n^{\IZ^d}\bigl(E\IZ^d;\bfL^{\langle -\infty \rangle}_R\bigr) 
& \cong &
\bigoplus_{i = 0}^d L_{n-i}^{\langle -\infty \rangle}(R)^{\binom{d}{i}}.
\label{HZd(EZd;bfL_R)}
\end{eqnarray}

Theorem~\ref{the:FJC_hyperbolic_virtually_Zd} implies that
$$\Wh_n(\IZ^d;R) := H^{\IZ^d}_n\bigl(E\IZ^d  \to \pt;\bfK_R\bigr)
\cong 
H^{\IZ^d}_n\bigl(E\IZ^d  \to \edub{\IZ^d};\bfK_R\bigr).$$
From~\cite[Corollary~2.10]{Lueck-Weiermann(2012)} it follows that
$$\bigoplus_{C \in \calmicyc} H^{\IZ^d}_n\bigl(E\IZ^d  \to E\IZ^d/C;\bfK_R\bigr)
\cong 
H^{\IZ^d}_n\bigl(E\IZ^d  \to \edub{\IZ^d};\bfK_R\bigr).
$$
Since $C \subseteq \IZ^d$ is maximal infinite cyclic,
$\IZ^d \cong C \oplus \IZ^{d-1}$.  Therefore the induction structure
and~\eqref{H_nZ(EZ_to_pt;bfK_R)_cong_NK_n(R)_oplus_NK_n(R)} imply
\begin{eqnarray}
H^{\IZ^d}_n\bigl(E\IZ^d  \to E(\IZ^d/C);\bfK_R\bigr)
& \cong &
H^{C \oplus \IZ^{d-1}}_n\bigl(EC \times E\IZ^{d-1}  \to E\IZ^{d-1};\bfK_R\bigr)
\nonumber \\
& \cong &
H_n^{C}\bigl(EC\times B\IZ^{d-1} \to B\IZ^{d-1};\bfK_R\bigr)
\nonumber \\
& \cong &
H^{C}_n\bigl((EC \to \pt) \times B\IZ^{d-1};\bfK_R\bigr)
\nonumber \\
& \cong &
\bigoplus_{i = 0}^{d-1} 
H^{C}_{n-i}\bigl(EC \to \pt;\bfK_R\bigr)^{\binom{d-1}{i}}
\nonumber \\
& \cong &
\bigoplus_{i = 0}^{d-1} 
\bigl(\NK_{n-i}(R) \oplus \NK_{n-i}(R)\bigr)^{\binom{d-1}{i}}.
\label{Zdrelativeterm}
\end{eqnarray}

Hence,
\begin{eqnarray*}
\Wh_n(\IZ^d;R) 
& \cong &
\bigoplus_{C \in \calmicyc}
\bigoplus_{i = 0}^{d-1} 
\bigl(\NK_{n-i}(R) \oplus \NK_{n-i}(R)\bigr)^{\binom{d-1}{i}}.
\end{eqnarray*}
From Theorem~\ref{the:FJC_hyperbolic_virtually_Zd} and
Lemma~\ref{lem:splitting_eub_deub}~\ref{lem:splitting_eub_deub:K}
\[ K_n(R[\IZ^d]) \; \cong \; H_n^{\IZ^d}(E\IZ^d;\bfK_R) \oplus
\Wh_n(\IZ^d;R), \] which, by~\eqref{HZd(EZd;bfk_R)}, is isomorphic to
\[ \left(\bigoplus_{i = 0}^{d} K_{n-i}(R)^{\binom{d}{i}} \right) \oplus
\Wh_n(\IZ^d;R). \] Finally,
\[ L_n^{\langle -\infty \rangle}(R[\IZ^d]) \; \cong \;
H_n^{\IZ^d}\bigl(E\IZ^d;\bfL^{\langle -\infty \rangle}_R\bigr) \; \cong \;
\bigoplus_{i = 0}^d L_{n-i}^{\langle -\infty \rangle}(R)^{\binom{d}{i}}\] by
Theorem~\ref{the:FJC_hyperbolic_virtually_Zd},
Remark~\ref{rem:role_of_type_I_and_II} and~\eqref{HZd(EZd;bfL_R)}.
\end{proof}


\subsection{Virtually finitely generated abelian groups}
\label{subsec:proof_Virtually_finitely_generated_abelian_groups}

For the remainder of the paper
we will use the notation introduced in 
Subsection~\ref{subsec:K-theory_in_the_case_of_a_free_conjugation_action}
and the following notation. For $C \in I$ there is an obvious extension 
\begin{eqnarray*}
& 1 \to C \to N_GC \xrightarrow{p_C} W_GC \to 1, &
\end{eqnarray*}
where $p_C$ is the canonical projection. We also have the extension
\begin{eqnarray}
&1 \to A/C \to W_GC \xrightarrow{\overline{q_C}} Q_C \to 1, &
\label{extension_A/C_W_GC_Q_C}
\end{eqnarray} 
which is induced by the given extension $1 \to A  \to G \xrightarrow{q} Q \to 1$.
Since $C \subseteq A$ is a maximal infinite cyclic subgroup, $A/C \cong \IZ^{d-1}$.

Notice that any infinite cyclic subgroup $C$ of $A$ is contained in a unique maximal
infinite cyclic subgroup $C_{\max}$ of $A$.  In particular for two maximal
infinite cyclic subgroups $C,D \subseteq A$, either $C \cap D =
\{0\}$ or $C = D$.  
Let
$$N_G[C] := \{g \in G \mid |gCg^{-1} \cap C| = \infty\}.$$
For every $C \in I$
$$N_G[C] = N_GC = q^{-1}(Q_C).$$
Consider the following equivalence relation on the set of infinite virtually
cyclic subgroups of $G$. We call $V_1$ and $V_2$ equivalent if and only
if $(A\cap V_1)_{\max} = (A\cap V_2)_{\max}$. Then for every infinite virtually
cyclic subgroup $V$ of $G$ there is precisely one $C \in I$
such that $V$ is equivalent to $gCg^{-1}$, for some $g \in G$.
We obtain from~\cite[Theorem~2.3]{Lueck-Weiermann(2012)} isomorphisms
\begin{eqnarray}
\hspace{-8mm} \bigoplus_{C \in I} H_n^{N_GC}\bigl(\eub{N_GC} \to p_C^*\eub{W_GC};\bfK_R\bigr)
& \cong &
H_n^G\bigl(\eub{G} \to \edub{G};\bfK_R\bigr);
\label{H_nG(eubG_to_edubG;bfK_R)_virtually_Zd}
\\
\hspace{-8mm} \bigoplus_{C \in I} H_n^{N_GC}\bigl(\eub{N_GC} \to p_C^*\eub{W_GC};
\bfL^{\langle - \infty\rangle}_R\bigr)
& \cong &
H_n^G\bigl(\eub{G} \to \edub{G};\bfL^{\langle - \infty\rangle}_R\bigr),
\label{H_nG(eubG_to_edubG;bfL-infty_R)_virtually_Zd}
\end{eqnarray}
where $p_C \colon N_GC = q^{-1}(Q_C) \to W_GC = q^{-1}(Q_C)/C$ is the canonical projection.

\begin{lemma}\label{lem:useful_G-homeo}
Let $f \colon G_1 \to G_2$ be a surjective group homomorphism.
Consider a subgroup $H\subset G_2$. Let $Y$ be a $G_1$-space and $Z$ be 
an $H$-space. Denote by  $f_H \colon f^{-1}(H) \to H$ the map induced by $f$.

Then there is a natural $G_1$-homeomorphism
$$G_1 \times_{f^{-1}(H)} \bigl(\res_{G_1}^{f^{-1}(H)} Y \times f_H^* Z\bigr) 
\xrightarrow{\cong} Y \times f^*(G_2 \times_H Z),$$
where $f_H^*$, $f^*$ and $\res_{G_1}^{f^{-1}(H)}$ denote restriction and the actions on 
products are the diagonal actions.
\end{lemma}
\begin{proof}
The map  sends $\bigl(g,(y,z)\bigr)$ to $\bigl(gy,(f(g),z)\bigr)$.
Its inverse sends $\bigl(y,(k,z)\bigr)$ to $\bigl(h,(h^{-1}y,z)\bigr)$
for any $h \in G_1$ with $f(h) = k$.
\end{proof}


\subsubsection{$K$-theory in the case of a free conjugation action}
\label{subsubsec:proof_K-theory_in_the_case_of_a_free_conjugation_action}

\begin{proof}[Proof of Theorem~\ref{the:K-theory_of_G_with_Aq_is_0}]
We prove assertions (i) and (ii) simultaneously by a direct computation.

By Theorem~\ref{the:FJC_hyperbolic_virtually_Zd}
$H_n^G\bigl(\edub{G};\bfK_R\bigr) \cong H_n^G\bigl(\pt;\bfK_R\bigr)$. Thus,
\begin{eqnarray*}
\Wh_n(G;R) 
& := & H_n^G\bigl(EG\to \pt;\bfK_R\bigr)\\
& \cong & H_n^G\bigl(EG\to \edub{G};\bfK_R\bigr)\\
& \cong & H_n^G\bigl(EG\to \eub{G};\bfK_R\bigr) \oplus H_n^G\bigl(\eub{G}\to \edub{G};\bfK_R\bigr), 
\end{eqnarray*}
using~(\ref{relative_K_splitting}). From~\cite[Lemma~6.3]{Lueck-Stamm(2000)} 
and~\cite[Corollary~2.11]{Lueck-Weiermann(2012)}, there is a $G$-pushout
\begin{eqnarray}
&
\xymatrix{
     \coprod_{F\in J} G\times_{F}EF
     \ar[d]^{\coprod_{F\in J} p} \ar[r]^-i & EG \ar[d] \\
     \coprod_{F\in J} G/F \ar[r] & \eub{G}
 }
&
\label{pushout_for_eubG}
\end{eqnarray}
where $J$ is a complete system of representatives of maximal finite subgroups of
$G$. This produces an isomorphism
\begin{eqnarray*}
\bigoplus_{F \in J} \Wh_n(F;R):=
\bigoplus_{F \in J} H_n^F\bigl(EF \to \pt;\bfK_R\bigr) 
\xrightarrow{\cong} 
H_n^G\bigl(EG \to \eub{G};\bfK_R\bigr).
\end{eqnarray*}
By~\eqref{H_nG(eubG_to_edubG;bfK_R)_virtually_Zd},
\[ \bigoplus_{C \in I} H_n^{N_GC}\bigl(\eub{N_GC} \to p_C^*\eub{W_GC};\bfK_R\bigr) \; \cong \;
H_n^G\bigl(\eub{G} \to \edub{G};\bfK_R\bigr). \]

Fix $C$ in $I$. Since the conjugation action of $Q$ on $A$ is free away from
$0$, it induces an embedding of $Q_C$ into $\aut(C)$. Hence $Q_C$ is either
trivial or isomorphic to $\IZ/2$. Therefore, $I = I_1 \amalg I_2$, where $ I_1 =
\bigl\{C \in I \mid Q_C = \{1\}\bigr\}$ and $I_2 = \bigl\{C \in I \mid Q_C =
\IZ/2\bigr\}$.

From~\cite[Corollary~2.10]{Lueck-Weiermann(2012)} there is an isomorphism
\begin{multline*}
\bigoplus_{C \in \calmicyc(A)} H^{N_AC}_n\bigl(EN_AC \to p_C^*EW_AC;\bfK_R\bigr)\\
\cong H^A_n\bigl(EA  \to \edub{A};\bfK_R\bigr)=\Wh_n(A;R),
\end{multline*}
where $p_C \colon N_AC = A \to W_AC = A/C$ denotes the projection. The conjugation action of $Q$ on $A$ induces an action of $Q$ on $H^A_n\bigl(EA  \to \edub{A};\bfK_R\bigr)$. 
By the definition of the index set $I$ and the subgroup $Q_C \subseteq Q$, we obtain a $\IZ Q$-isomorphism
$$\bigoplus_{C \in I} \IZ Q \otimes_{\IZ[Q_C]} H^A_n\bigl(EA \to p_C^*E(A/C);\bfK_R\bigr)
\cong 
H^{A}_n\bigl(EA  \to \edub{A};\bfK_R\bigr),
$$
and hence an isomorphism
$$\bigoplus_{C \in I} \IZ  \otimes_{\IZ[Q_C]} H^A_n\bigl(EA \to p_C^*E(A/C);\bfK_R\bigr)
\cong 
\IZ  \otimes_{\IZ Q} H^{A}_n\bigl(EA  \to \edub{A};\bfK_R\bigr).
$$

If $Q_C = \{1\}$, then $N_GC = A\cong C\oplus \IZ^{d-1}$ and $W_GC\cong \IZ^{d-1}$. 
Thus, by~(\ref{Zdrelativeterm}),
\begin{eqnarray*}
H_n^{N_GC}\bigl(\eub{N_GC} \to p_C^*\eub{W_GC};\bfK_R\bigr)
& \cong & H^A_n\bigl(EA \to p_C^*E(A/C);\bfK_R\bigr)
\\
& \cong &
H_n^{C \oplus \IZ^{d-1}}\bigl(EC \times E\IZ^{d-1} \to E\IZ^{d-1};\bfK_R\bigr)
\\
& \cong &
\bigoplus_{i = 0}^{d-1} 
\bigl(\NK_{n-i}(R) \oplus \NK_{n-i}(R)\bigr)^{\binom{d-1}{i}}.
\end{eqnarray*}
Therefore, the proof of the theorem will be finished once it is established that
for $C\in I_2$, 
\begin{multline*}
H_n^{N_GC}\bigl(\eub{N_GC} \to p_C^*\eub{W_GC};\bfK_R\bigr)\\
\cong \IZ  \otimes_{\IZ[Q_C]} H^A_n\bigl(EA \to p_C^*E(A/C);\bfK_R\bigr)
\cong \bigoplus_{i = 0}^{d-1} \NK_{n-i}(R)^{\binom{d-1}{i}}.
\end{multline*}

Assume that $Q_C=\IZ/2$. Since the projection $\eub{N_GC} \times
p_C^*\eub{W_GC}\to \eub{N_GC}$ is an $N_GC$-homotopy equivalence,
$H_n^{N_GC}\bigl(\eub{N_GC} \to p_C^*\eub{W_GC};\bfK_R\bigr)$ is isomorphic to
$H_n^{N_GC}\bigl(\eub{N_GC} \times p_C^*\eub{W_GC} \to
p_C^*\eub{W_GC};\bfK_R\bigr)$, which, by
Lemma~\ref{lem:intermediate_isomorphism}(i) below, is isomorphic to
\[ H_n^{N_GC}\bigl(\eub{N_GC} \times p_C^*(EW_GC \times \overline{q_C}^*EQ_C)
\to p_C^*(EW_GC \times \overline{q_C}^*EQ_C) ;\bfK_R\bigr). \] 


Sending a $Q_C$-$CW$-complex $Y$ to the $\IZ$-graded abelian group
$$H_*^{N_GC}\bigl(\eub{N_GC} \times  p_C^*(EW_GC \times \overline{q_C}^*Y) 
\to p_C^*(EW_GC \times \overline{q_C}^*Y) ;\bfK_R\bigr)$$ yields a
$Q_C$-homology theory.  Since $EQ_C$ is a free $Q_C$-$CW$-complex, there is an
equivariant Atiyah-Hirzebruch spectral sequence converging to
$$H_{i+j}^{N_GC}\bigl(\eub{N_GC} \times  p_C^*(EW_GC \times \overline{q_C}^*EQ_C) 
\to p_C^*(EW_GC \times \overline{q_C}^*EQ_C) ;\bfK_R\bigr)
$$
whose $E^2$-term is given by
\begin{multline*}
E^2_{i,j} = H_i^{Q_C}\biggl(EQ_C; 
H_{j}^{N_GC}\bigl(\eub{N_GC} \times  p_C^*(EW_GC \times \overline{q_C}^*Q_C) 
\\
\to  p_C^*(EW_GC \times \overline{q_C}^*Q_C) ;\bfK_R\bigr)\bigg).
\end{multline*}

Lemma~\ref{lem:intermediate_isomorphism}~\ref{lem:intermediate_isomorphism:(ii)} implies that:
\begin{eqnarray*}
E^2_{i,j} 
& \cong &
H_i^{Q_C}\left(EQ_C;
H^A_j\bigl(EA \to p_C^*E(A/C);\bfK_R\bigr)\right)
\\
& \cong &
H_i^{Q_C}\left(EQ_C;
\IZ[Q_C]  \otimes_{\IZ}\left(\IZ  \otimes_{\IZ[Q_C]}H^A_j\bigl(EA \to p_C^*E(A/C);\bfK_R\bigr)\right)\right)
\\
& \cong &
H_i^{\{1\}}\bigl(\res_{Q_C}^{\{1\}} EQ_C;\IZ  \otimes_{\IZ[Q_C]}H^A_j\bigl(EA \to p_C^*E(A/C);\bfK_R\bigr)\bigr)
\\
& \cong &
\begin{cases}
\IZ  \otimes_{\IZ[Q_C]}H^A_j\bigl(EA \to p_C^*E(A/C);\bfK_R\bigr) & i = 0
\\
\{0\} & i \not= 0.
\end{cases}
\\
& \cong &
\begin{cases}
\bigoplus_{l = 0} ^{d-1}   \NK_{j-l}(R)^{\binom{d-1}{l}} & i = 0
\\
\{0\} & i \not= 0.
\end{cases}
\end{eqnarray*}
Hence, the Atiyah-Hirzebruch spectral sequence collapses at $E^2$. Therefore,
\begin{multline*}
H_n^{N_GC}\bigl(\eub{N_GC} \times  p_C^*(EW_GC \times \overline{q_C}^*EQ_C) 
\to p_C^*(EW_GC \times \overline{q_C}^*EQ_C) ;\bfK_R\bigr)
\\ \cong \IZ  \otimes_{\IZ[Q_C]} H^A_n\bigl(EA \to p_C^*E(A/C);\bfK_R\bigr)
\cong \bigoplus_{l = 0} ^{d-1}  \NK_{n-l}(R)^{\binom{d-1}{l}}.
\end{multline*}
This completes the proof of Theorem~\ref{the:K-theory_of_G_with_Aq_is_0} once
we have proved Lemma~\ref{lem:intermediate_isomorphism}.
\end{proof}

\begin{lemma}\label{lem:intermediate_isomorphism}
Let $G$ be a group satisfying the assumptions of
Theorem~\ref{the:K-theory_of_G_with_Aq_is_0}, and let $C$ be a maximal cyclic
subgroup of $G$ such that $Q_C = \IZ/2$. Then:
\begin{enumerate}
\item \label{lem:intermediate_isomorphism:(i)}
  For every $n \in \IZ$, the projection $EW_GC \times \overline{q_C}^*EQ_C
  \to \eub{W_GC}$ induces a bijection,
  \begin{multline*}
    H_n^{N_GC}\bigl(\eub{N_GC} \times p_C^*(EW_GC \times \overline{q_C}^*EQ_C)
    \to p_C^*(EW_GC \times \overline{q_C}^*EQ_C) ;\bfK_R\bigr)
    \\
    \xrightarrow{\cong} H_n^{N_GC}\bigl(\eub{N_GC} \times p_C^*\eub{W_GC} \to
    p_C^*\eub{W_GC};\bfK_R\bigr).
  \end{multline*}
\item \label{lem:intermediate_isomorphism:(ii)} There are isomorphisms of $\IZ[Q_C]$-modules
  \begin{multline*}
    H_{j}^{N_GC}\bigl(\eub{N_GC} \times p_C^*(EW_GC \times
    \overline{q_C}^*Q_C) \to p_C^*(EW_GC \times \overline{q_C}^*Q_C)
    ;\bfK_R\bigr) \\ \cong H^A_j\bigl(EA \to p_C^*E(A/C);\bfK_R\bigr)  \cong \bigoplus_{l = 0} ^{d-1} \bigl(\IZ[Q_C]
    \otimes_{\IZ} \NK_{j-l}(R)\bigr)^{\binom{d-1}{l}}.
  \end{multline*}
\end{enumerate}
\end{lemma}

\begin{proof}
(i) Note that the conjugation action of $Q_C=\IZ/2$ on $A/C$ is free away from
$0 \in A/C$. To see this, let $q$ be a nontrivial element in $Q$, and let $aC
\in A/C$ such that $aC$ is fixed under the conjugation action with $q$. Then
there is a $c \in C$ such that $\rho(q)(a) = a + c$, where $\rho \colon Q \to
\aut(A)$ is given by the conjugation action of $Q$ in $A$. This implies
$\rho(q)(2a) = 2a + 2c$. Thus,
$$\rho(q)(2a + c) = \rho(q)(2a) + \rho(q)(c) = \rho(q)(2a) - c = 2a +2c - c = 2a+c.$$
By assumption, the $Q$-action on $A$ is free away from $0 \in A$, so $2a+c = 0$
and thus $2a \in C$. Since $C$ is maximal infinite cyclic, $aC = 0$ in
$A/C$. Hence, the conjugation action of $Q_C$ on $A/C$ is free away from $0 \in A/C$.

Since $(A/C)^{Q_C} = \{0\}$, the extension~\eqref{extension_A/C_W_GC_Q_C} has a
section, and $W_GC$ is isomorphic to the semidirect product $\IZ^{d-1} \rtimes
\IZ/2$ with respect to the involution $-\id \colon \IZ^{d-1} \to \IZ^{d-1}$.
Let $J_C$ be a complete system of representatives of the conjugacy classes of
maximal finite subgroups of $W_GC$.  Then every element in $J_C$ is isomorphic
to $\IZ/2$ and $J_C$ contains $2^{d-1}$ elements (see
Remark~\ref{rem:cardinality_of_J}).  From~\cite[Lemma~6.3]{Lueck-Stamm(2000)}
and~\cite[Corollary~2.11]{Lueck-Weiermann(2012)}, there is a $W_GC$-pushout
\begin{eqnarray}
&
\xymatrix{
  \coprod_{H\in J_C} W_GC\times_{H}EH
  \ar[d]^{\coprod_{H\in J_C} p} \ar[r]^-i & EW_GC \ar[d] \\
  \coprod_{H\in J_C} W_GC/H \ar[r] & \eub{W_GC}
}
&
\label{pushout_for_eub(W_GC)}
\end{eqnarray}

Let $H \in J_C$.  Given any $N_GC$-space $Y$ and any $H$-space $Z$,
Lemma~\ref{lem:useful_G-homeo} produces an $N_GC$-homeomorphism
\begin{eqnarray}
N_GC \times_{p_C^{-1}(H)} \bigl(\res_{N_GC}^{p_C^{-1}(H)}(Y) 
\times p_C^* Z\bigr) 
\xrightarrow{\cong}  Y \times p_C^*(W_GC\times_{H}Z),
\label{NGC_homeo}
\end{eqnarray} 
where $p_C^*$ denotes restriction for both $p_C \colon N_GC \to W_GC$ and
the homomorphism $p^{-1}_C(H) \to H$ induced by $p_C$. Since
$\res_{N_GC}^{p_C^{-1}(H)}(\eub{N_GC})$ is a $p_C^{-1}(H)$-$CW$-model for
$\eub{p_C^{-1}(H)}$, we obtain identifications of $N_GC$-spaces
\begin{eqnarray*}
\eub{N_GC} \times p_C^*(W_GC\times_{H} EH) 
& \cong &
N_GC \times_{p_C^{-1}(H)} \bigl(\eub{p_C^{-1}(H)} \times p_C^*EH\bigr);
\\
\eub{N_GC} \times p_C^*(W_GC/H) 
& \cong &
N_GC \times_{p_C^{-1}(H)} \bigl(\eub{p_C^{-1}(H)}\bigr);
\\
p_C^*(W_GC\times_{H} EH) 
& \cong &
N_GC \times_{p_C^{-1}(H)}  p_C^*EH;
\\
p_C^*(W_GC/H) 
& \cong &
N_GC \times_{p_C^{-1}(H)} \pt,
\end{eqnarray*}
by substituting $Y=\eub{N_GC}$ or $Y=\pt$, and $Z=EH$ or $Z=\pt$
into~\eqref{NGC_homeo}.  Using the induction structure of the equivariant
homology theory $H_*^?(-;\bfK_R)$ and the Five-Lemma, these identifications
imply
\begin{multline}
H_n^{N_GC}\bigl(\eub{N_GC} \times p_C^*(W_GC\times_{H}EH) \to
p_C^*(W_GC\times_{H}EH);\bfK_R\bigr)
\\
\xrightarrow{\cong} H_n^{p_C^{-1}(H)}\bigl(\eub{p_C^{-1}(H)} \times p_C^*EH
\to p_C^*EH;\bfK_R\bigr)
\label{map1_for_K_and_W_GC}
\end{multline}
and
\begin{multline}
H_n^{p_C^{-1}(H)}\bigl(\eub{p_C^{-1}(H)} \to \pt;\bfK_R\bigr)
\\
\xleftarrow{\cong} H_n^{N_GC}\bigl(\eub{N_GC} \times p_C^*(W_GC/H) \to
p_C^*(W_GC/H);\bfK_R\bigr)
\label{map2_for_K_and_W_GC}
\end{multline}
are bijective for every $n \in \IZ$.

Consider the following $p_C^{-1}(H)$-pushout, where the left
vertical arrow is an inclusion of $p_C^{-1}(H)$-$CW$-complexes and the upper horizontal arrow is cellular.
$$\xymatrix{Ep_C^{-1}(H)        
     \ar[d] \ar[r] & p_C^*EH \ar[d] \\
      \eub{p_C^{-1}(H)}  \ar[r] & X
 }
$$
Then $X$ is a $p_C^{-1}(H)$-$CW$-complex, and for every subgroup
$K \subseteq p_C^{-1}(H)$, the above $p_C^{-1}(H)$-pushout induces
a pushout of $CW$-complexes:
$$\xymatrix{Ep_C^{-1}(H)^K        
     \ar[d] \ar[r] & (p_C^*EH)^K \ar[d] \\
      \eub{p_C^{-1}(H)}^K  \ar[r] & X^K
 }
$$
If $K = \{1\}$, then the spaces $Ep_C^{-1}(H)^K$, $(p_C^*EH)^K$ and
$\eub{p_C^{-1}(H)}^K$ are contractible, and thus $X^K$ is contractible.  If $K$
is a non-trivial finite subgroup, then the spaces $Ep_C^{-1}(H)^K$ and
$(p_C^*EH)^K$ are empty, and the space $\eub{p_C^{-1}(H)}^K$ is
contractible. Hence $X^K$ is contractible. If $K$ is an infinite subgroup of
$C$, then the spaces $Ep_C^{-1}(H)^K$ and $\eub{p_C^{-1}(H)}^K$ are empty, and
the space $(p_C^*EH)^K$ is contractible since $p_C(K)=\{1\}$. Therefore $X^K$ is
contractible. If $K$ is infinite and not contained in $C$, then the spaces
$Ep_C^{-1}(H)^K$, $(p_C^*EH)^K$ and $\eub{p_C^{-1}(H)}^K$ are all empty, and so
$X^K$ is also empty.  Since $p_C^{-1}(H) \cong D_{\infty}$, an infinite
virtually cyclic subgroup $K$ of $p_C^{-1}(H)$ is of type I if and only if it is
an infinite subgroup of $C$. Hence, $X$ is a model for
$\EGF{{p_C^{-1}(H)}}{\calvcyc_I}$. Notice that $\eub{p_C^{-1}(H)} \times
p_C^*EH$ is a model for $Ep_C^{-1}(H)$. Therefore, we have the following
cellular $p_C^{-1}(H)$-pushout, where the left vertical arrow is an inclusion of
$p_C^{-1}(H)$-$CW$-complexes and the upper horizontal arrow is cellular.
$$\xymatrix{\eub{p_C^{-1}(H)} \times p_C^*EH       
\ar[d] \ar[r] & p_C^*EH \ar[d] \\
\eub{p_C^{-1}(H)} \ar[r] & \EGF{{p_C^{-1}(H)}}{\calvcyc_I} }
$$
This induces an isomorphism
\begin{multline}
H_n^{p_C^{-1}(H)}\bigl(\eub{p_C^{-1}(H)} \times p_C^*EH \to
p_C^*EH;\bfK_R\bigr)
\\
\xrightarrow{\cong} H_n^{p_C^{-1}(H)}\bigl(\eub{p_C^{-1}(H)} \to
\EGF{{p_C^{-1}(H)}}{\calvcyc_I};\bfK_R\bigr),
\label{map_for_K_and_W_GC_overPC-1(H)_vcycI}
\end{multline}
for every $n\in\IZ$.  Since $p_C^{-1}(H) \cong D_{\infty}$ is virtually cyclic,
Remark~\ref{rem:role_of_type_I_and_II} implies that the map
\begin{eqnarray*}
H_n^{p_C^{-1}(H)}\bigl(\EGF{p_C^{-1}(H)}{\calvcyc_I};\bfK_R\bigr)
&\xrightarrow{\cong} &
H_n^{p_C^{-1}(H)}\bigl(\pt;\bfK_R\bigr)
\end{eqnarray*}
is bijective for every $n \in \IZ$, and so by
\eqref{map_for_K_and_W_GC_overPC-1(H)_vcycI},
\begin{multline}
H_n^{p_C^{-1}(H)}\bigl(\eub{p_C^{-1}(H)} \times p_C^*EH \to
p_C^*EH;\bfK_R\bigr)
\\
\xrightarrow{\cong} H_n^{p_C^{-1}(H)}\bigl(\eub{p_C^{-1}(H)} \to
\pt;\bfK_R\bigr)
\label{map_for_K_and_W_GC_overPC-1(H)}
\end{multline}
is bijective for every $n \in
\IZ$. Therefore,~\eqref{map1_for_K_and_W_GC},~\eqref{map2_for_K_and_W_GC}
and~\eqref{map_for_K_and_W_GC_overPC-1(H)} imply that
\begin{multline}
H_n^{N_GC}\bigl(\eub{N_GC} \times p_C^*(W_GC\times_{H}EH) \to
p_C^*(W_GC\times_{H}EH);\bfK_R\bigr)
\\
\xrightarrow{\cong} H_n^{N_GC}\bigl(\eub{N_GC} \times p_C^*(W_GC/H) \to
p_C^*(W_GC/H);\bfK_R\bigr)
\label{map_for_K_and_W_GC}
\end{multline}
is an isomorphism for every $n\in\IZ$.

We obtain a $W_GC$-homology theory by assigning to a $W_GC$-$CW$-complex $Z$ the
$\IZ$-graded abelian group  
$H_n^{N_GC}\bigl(\eub{N_GC} \times  p_C^*Z\to p_C^*Z;\bfK_R\bigr)$.
From the Mayer-Vietoris sequence associated to the 
$W_GC$-pushout~\eqref{pushout_for_eub(W_GC)} and the 
bijectivity of the map~\eqref{map_for_K_and_W_GC}, there is an isomorphism
\begin{multline}
H_n^{N_GC}\bigl(\eub{N_GC} \times  p_C^*EW_GC \to p_C^*EW_GC;\bfK_R\bigr)
\\
\xrightarrow{\cong} 
H_n^{N_GC}\bigl(\eub{N_GC} \times  p_C^*\eub{W_GC} \to p_C^*\eub{W_GC};\bfK_R\bigr).
\label{iso_for_H_nN_GC(eub(N_GC)_times_p_Ceub(W_GC)_to_p_Ceub(W_GC);bfK_R)}
\end{multline}
Since the projection $EW_GC \times \overline{q_C}^*EQ_C \to EW_GC $ is a
$W_GC$-homotopy equivalence, the desired isomorphism now follows
from \eqref{iso_for_H_nN_GC(eub(N_GC)_times_p_Ceub(W_GC)_to_p_Ceub(W_GC);bfK_R)}.
\vspace{.25cm}

(ii) By Lemma~\ref{lem:useful_G-homeo}
\begin{eqnarray}
\lefteqn{H_{j}^{N_GC}\bigl(\eub{N_GC} \times  p_C^*(EW_GC \times \overline{q_C}^*Q_C) 
\to p_C^*(EW_GC \times \overline{q_C}^*Q_C) ;\bfK_R\bigr)}
& &
\nonumber
\\
& \cong &
H_{j}^{N_GC}\bigl(\eub{N_GC} \times  p_C^*(W_GC \times_{A/C} \res_{W_GC}^{A/C}EW_GC) 
\nonumber
\\
& & \qquad\to p_C^*(W_GC \times_{A/C} \res_{W_GC}^{A/C}EW_GC);\bfK_R\bigr)
\nonumber
\\
& \cong &
H_{j}^{N_GC}\biggl(N_GC \times_{A} \bigl(\res_{N_GC}^{A}\eub{N_GC} \times 
\overline{p}_C^*\res_{W_GC}^{A/C}EW_GC)\bigr)
\nonumber
\\ & & 
\qquad \to N_GC \times_{A} \bigl(\overline{p}_C^*\res_{W_GC}^{A/C}EW_GC\bigr);\bfK_R\bigg).
\label{zero-th-identification}
\end{eqnarray}
The generator $t$ of $Q_C = \IZ/2$ acts on 
$$H_{j}^{N_GC}\bigl(\eub{N_GC} \times  p_C^*(EW_GC \times \overline{q_C}^*Q_C) 
\to p_C^*(EW_GC \times \overline{q_C}^*Q_C) ;\bfK_R\bigr)$$
by the square of $N_GC$-maps
$$\comsquare{\eub{N_GC} \times  p_C^*(EW_GC \times \overline{q_C}^*Q_C)}
{}
{p_C^*(EW_GC \times \overline{q_C}^*Q_C)}
{\id \times p_C^*(\id \times r_t)}{p_C^*(\id \times r_t)}
{\eub{N_GC} \times  p_C^*(EW_GC \times \overline{q_C}^*Q_C)}
{}
{p_C^*(EW_GC \times \overline{q_C}^*Q_C)}
$$
where $r_t \colon Q_C \to Q_C$ is right multiplication with $t$.

To simplify notation, let $X=\res_{N_GC}^{A}\eub{N_GC}$ and $Y=\overline{p}_C^*\res_{W_GC}^{A/C}EW_GC$. Choose an element $\gamma \in N_GC$ that is mapped to $t$ by $q_C \colon N_GC \to Q_C$, and equip 
$$H_{j}^{N_GC}\bigl(N_GC \times_{A} (X \times Y) \to N_GC \times_{A} Y;\bfK_R\bigr)$$
with the $Q_C$-action coming from square of $N_GC$-maps
$$
\comsquare{N_GC \times_{A} \bigl(X \times Y\bigr)}
{}
{N_GC \times_{A} Y}
{r_{\gamma} \times l_{\gamma^{-1}}}
{r_{\gamma} \times l_{\gamma^{-1}}}
{N_GC \times_{A} \bigl(X \times Y\bigr)}
{}
{N_GC \times_{A} Y}
$$
where $l_{\gamma^{-1}}$ is given by $(x,y) \mapsto (\gamma^{-1} \cdot x, p_C(\gamma^{-1}) \cdot y)$ 
for $x \in X$ and $y \in Y$ and $r_{\gamma}$ is right
multiplication by $\gamma$. It is straightforward to check that the 
isomorphism~\eqref{zero-th-identification} is compatible with this $Q_C$-action.

Recall that $Q_C$ acts freely on $A$ away from the identity. This implies that if $t \in Q_C$ is the generator and $a \in A$, then $t \cdot (a + ta) = a + ta$, and so $a + ta = 0$. Therefore the $Q_C$-action on $A$ is given by $-\id_A \colon
A \to A$. Since the action of $t$ on $A$ is defined by conjugation by
$\gamma^{-1}$, the following diagram commutes.
\[
\xymatrix{
	N_GC
		\ar[r]^{c(\gamma^{-1})} &
	N_GC
		\\
	A
		\ar[u] \ar[r]^{-\id_A} &
	A
		\ar[u] 
}
\]
For any $A$-CW complex $Z$, this yields the commutative square
\[
\xymatrix@C=6em{
	H_{j}^{N_GC}\bigl(N_GC \times_{A} Z;\bfK_R\big)
		\ar[r]^-{f_*\circ \ind_{c(\gamma^{-1})}}_-\cong &
	H_{j}^{N_GC}\bigl(N_GC \times_{A} (A\times_{-\id_A}Z);\bfK_R\big)
		\\
	H_{j}^A(Z;\bfK_R)
		\ar[u]^{\ind_A^{N_GC}}_-\cong \ar[r]^-{\ind_{-\id_A}}_-\cong &
	H_{j}^A(A\times_{-\id_A}Z;\bfK_R),
		\ar[u]^{\ind_A^{N_GC}}_-\cong
}
\]
where $\ind_A^{N_GC}$ denotes induction with respect to inclusion, and $f_*$ is induced by the map $f:N_GC\times_{c(\gamma^{-1})}(N_GC \times_{A} Z) \to N_GC \times_{A} (A\times_{-\id_A}Z)$, which sends $\big(k,(g,z\big)$ to $\big(kc(\gamma^{-1})(g),(e,z)\big)$. From the axioms of an induction structure, we also have that $\ind_{c(\gamma^{-1})}$ is induced by the map $f_2: N_GC \times_{A} Z \to N_GC\times_{c(\gamma^{-1})}(N_GC \times_{A} Z) $, which sends $(k,z)$ to $\big(e,(\gamma k,z)\big)$ (see~\cite[Section~1]{Lueck(2002b)}).

When $Z=X\times Y$ or $Z=Y$, we can define a map $l:A\times_{-\id_A}Z \to Z$ such that $l(a,z)=a\gamma^{-1}z$. Since $\ind_A^{N_GC}$ is natural, this produces the commutative square
\[
\xymatrix@C=6em{
	H_{j}^{N_GC}\bigl(N_GC \times_{A} (A\times_{-\id_A}Z);\bfK_R\big)
		\ar[r]^-{l'_*}_-\cong &
	H_{j}^{N_GC}\bigl(N_GC \times_{A} Z;\bfK_R\big)
		\\
	H_{j}^A(A\times_{-\id_A}Z;\bfK_R)
		\ar[u]^{\ind_A^{N_GC}}_-\cong \ar[r]^-{l_*}_-\cong &
	H_{j}^A(Z;\bfK_R),
		\ar[u]^{\ind_A^{N_GC}}_-\cong
}
\]
where $l':N_GC \times_{A} (A\times_{-\id_A}Z) \to N_GC \times_{A} Z$ maps $\big(k (a,z)\big)$ to $(k,a\gamma^{-1}z)$. Notice that the composition $l'_*\circ f_* \circ \ind_{c(\gamma^{-1})}=l'_*\circ f_* \circ (f_2)_*$ is precisely the map $r_{\gamma}\times l_{\gamma^{-1}}$ used to define the $Q_C$-action on $H_{j}^{N_GC}\bigl(N_GC \times_{A} (X \times Y) \to N_GC \times_{A} Y;\bfK_R\bigr)$. Thus, the equivalence
\begin{eqnarray*} 
H_{j}^{N_GC}\bigl(N_GC \times_{A} (X \times Y) \to N_GC \times_{A} Y;\bfK_R\big) &  \xleftarrow{\ind_A^{N_GC}} &
H_{j}^{A}\bigl(X \times Y \to Y;\bfK_R\bigr) 
\end{eqnarray*}
is compatible with the $Q_C$-action, where the action on $H_{j}^{A}\bigl(X \times Y \to Y;\bfK_R\bigr)$ is induced by the composition $l_*\circ \ind_{-\id_A}$.

Notice that $X$ is a model for $EA$ and $Y$ is a model for $ED$. Therefore, the map $l$ defined above is unique up to $A$-homotopy. Since $C \subseteq A$ is maximal cyclic, there is a $D \subseteq A$ such that $C \oplus D \cong A$ and $D \cong \IZ^{d-1}$. Thus, the identification
\begin{eqnarray*}
H_{j}^{A}\bigl(X \times Y \to Y;\bfK_R\bigr)
& \cong &
H_{j}^{C \oplus D}\bigl((EC \times ED) \times ED \to ED;\bfK_R\bigr)
\\
& \cong &
H_{j}^{C \oplus D}\big(EC \times ED \to ED;\bfK_R\bigr)
\end{eqnarray*}
agrees with the $Q_C$-action since the action on $H_{j}^{C \oplus D}\big(EC \times ED \to ED;\bfK_R\bigr)$ is induced by $\big(l_C,l_D)_*\circ \ind_{(-\id_C\oplus -\id_D)}$, where $l_C:\ind_{-\id_C}EC \to EC$ denotes the unique $C$-equivariant map (up to $C$-homotopy) and $l_D:\ind_{-\id_D}ED \to ED$ denotes the unique $D$-equivariant map (up to $D$-homotopy). Note that this establishes the first desired isomorphism of $\IZ[Q_C]$-modules:
\begin{multline*}
 H_{j}^{N_GC}\bigl(\eub{N_GC} \times p_C^*(EW_GC \times
    \overline{q_C}^*Q_C) \to p_C^*(EW_GC \times \overline{q_C}^*Q_C)
    ;\bfK_R\bigr) \\ \cong H^A_j\bigl(EA \to p_C^*E(A/C);\bfK_R\bigr).
\end{multline*}


Using the freeness of the action of $D$,
\begin{eqnarray*}
H_{j}^{C \oplus D}\big(EC \times ED \to ED;\bfK_R\bigr)
& \cong &
H_{j}^{C}\big(EC\times BD \to BD;\bfK_R\bigr).
\end{eqnarray*}
As in~\eqref{Zdrelativeterm}, since $BD$ is the $(d-1)$-dimensional torus,
\begin{eqnarray*}
H_j^{C}\bigl(BD \times EC \to BD;\bfK_R\bigr)
& \cong &
\bigoplus_{l = 0} ^{d-1} H_{j-l}^{C}\bigl(EC \to \pt;\bfK_R\bigr)^{\binom{d-1}{l}} 
\\
& \cong &
\bigoplus_{l = 0} ^{d-1}\bigl(\NK_{j-l}(R) \oplus \NK_{j-l}(R)\bigr)^{\binom{d-1}{l}}.
\end{eqnarray*}
By Remark~\ref{K-theory_induction}, the $Q_C$-action on $H_*^{C}\bigl(EC \to \pt;\bfK_R\bigr)$, induced by $(l_C)_*\circ \ind_{-\id_C}$, 
corresponds to interchanging the two copies of $NK_*(R)$ in the decomposition
$$ 
H_*^{C}\bigl(EC \to \pt;\bfK_R\bigr)
\xrightarrow{\cong} 
\NK_*(R) \oplus \NK_*(R).
$$
Thus, we obtain an isomorphism
of $\IZ[Q_C]$-modules
\begin{eqnarray*}
\bigoplus_{l = 0} ^{d-1} \bigl(\IZ[Q_C] \otimes_{\IZ} \NK_{j-l}(R)\bigr)^{\binom{d-1}{l}} 
& \xrightarrow{\cong} &
H_j^{C}\bigl(BD \times EC \to BD;\bfK_R\bigr),
\end{eqnarray*}
which yields the desired result.
\end{proof}

\begin{proof}[Proof of Theorem~\ref{the:K-theory_for_regular_R_and_R_is_Z}]
\ \\~\ref{the:K-theory_for_regular_R_and_R_is_Z:R_regular} Since $R$ is a
regular ring, $\NK_j(R) = \{0\}$ for every $j \in \IZ$. Therefore,
$$\bigoplus_{F \in J} \Wh_n(F;R) \cong  \Wh_n(G;R)$$
for every $n \in \IZ$ by
Theorem~\ref{the:K-theory_of_G_with_Aq_is_0}. Furthermore, $K_n(R) = 0$ for
every $n \le -1$.  Thus, for any group $\Gamma$, the Atiyah-Hirzebruch spectral
sequence implies that
$$H_n(B\Gamma;\bfK(R)) = 0 \quad \text{for} \; n \le -1,$$
and the edge homomorphism induces an isomorphism
$$H_0(B\Gamma;\bfK(R)) \xrightarrow{\cong}  K_0(R).$$
Hence,
$$\Wh_n(\Gamma;R) \cong K_n(R\Gamma) \quad \text{for} \; n \le -1$$
and
$$\Wh_0(\Gamma;R) \cong \coker\bigl(K_0(R) \to K_0(R\Gamma)\bigr).$$%
\ref{the:K-theory_for_regular_R_and_R_is_Z:R_Dedekind} By
Carter~\cite{Carter(1980)},
$$K_n(RF) \cong \{0\} \quad \text{for} \; n \le -2$$
for every $F \in J$. Now apply
assertion~\ref{the:K-theory_for_regular_R_and_R_is_Z:R_regular} to complete the
proof.
\end{proof}


\subsubsection{$L$-theory in the case of a free conjugation action}
\label{subsubsec:proof_L-theory_in_the_case_of_a_free_conjugation_action}

\begin{proof}[Proof of Theorem~\ref{the:L-theory_of_G_with_Aq_is_0}]

As in the $K$-theory proof, Theorem~\ref{the:FJC_hyperbolic_virtually_Zd} and~(\ref{relative_L_splitting}) imply:
\begin{eqnarray*}
\cals_n^{\per,\langle -\infty \rangle}(G;R)
& := & H_n^G\bigl(EG\to \pt;\bfL^{\langle -\infty \rangle}_R\bigr)\\
& \cong & H_n^G\bigl(EG\to \edub{G};\bfL^{\langle -\infty \rangle}_R\bigr)\\
& \cong & H_n^G\bigl(EG\to \eub{G};\bfL^{\langle -\infty \rangle}_R\bigr) 
\oplus H_n^G\bigl(\eub{G}\to \edub{G};\bfL^{\langle -\infty \rangle}_R\bigr).
\end{eqnarray*}
From the $G$-pushout~(\ref{pushout_for_eubG}),
\begin{eqnarray*}
\bigoplus_{F \in J} \cals_n^{\per,\langle -\infty \rangle}(F;R):=
\bigoplus_{F \in J} H_n^F\bigl(EF \to \pt;\bfL^{\langle -\infty \rangle}_R\bigr) 
\cong
H_n^G\bigl(EG \to \eub{G};\bfL^{\langle -\infty \rangle}_R\bigr).
\end{eqnarray*}
By~\eqref{H_nG(eubG_to_edubG;bfL-infty_R)_virtually_Zd},
\[ \bigoplus_{C \in I} H_n^{N_GC}\bigl(\eub{N_GC} \to p_C^*\eub{W_GC};\bfL^{\langle -\infty \rangle}_R\bigr) \; \cong \;
H_n^G\bigl(\eub{G} \to \edub{G};\bfL^{\langle -\infty \rangle}_R\bigr). \]
Recall from the proof of Theorem~\ref{the:K-theory_of_G_with_Aq_is_0} that $I = I_1 \amalg I_2$. 

If $C\in I_1$, then $N_GC = A\cong C\oplus \IZ^{d-1}$ and $W_GC\cong \IZ^{d-1}$. 
Therefore, as in~(\ref{Zdrelativeterm}), 
$H_n^{N_GC}\bigl(\eub{N_GC} \to p_C^*\eub{W_GC};\bfL^{\langle -\infty \rangle}_R\bigr)$ is isomorphic to:
\[
H_n^{C \oplus \IZ^{d-1}}\bigl(EC \times E\IZ^{d-1} \to E\IZ^{d-1};\bfL^{\langle -\infty \rangle}_R\bigr)
\; \cong \;
\bigoplus_{i = 0}^{d-1} 
H_{n-i}^{C}\bigl(EC  \to \pt;\bfL^{\langle -\infty \rangle}_R\bigr)\;=\;0
\]
by Remark~\ref{rem:role_of_type_I_and_II}. 

Now assume $C\in I_2$. We obtain a $W_GC$-homology theory
by assigning to a $W_GC$-$CW$-complex $Z$ the $\IZ$-graded abelian group
$$
H_*^{N_GC}\bigl(\eub{N_GC} \times p_C^*Z \to  p_C^*Z;\bfL_R^{\langle - \infty \rangle}\bigr).
$$
Consider any free $W_GC$-$CW$-complex $Y$. Then there is an equivariant
Atiyah-Hirzebruch spectral sequence converging to 
$$
H_{i+j}^{N_GC}\bigl(\eub{N_GC} \times p_C^*Y 
\to  p_C^*Y;\bfL_R^{\langle - \infty \rangle}\bigr),
$$
whose $E^2$-term is given by
$$
E^2_{i,j} = H_i^{W_GC}\left(Y;H_{j}^{N_GC}\bigl(
\eub{N_GC} \times p_C^*W_GC \to  p_C^*W_GC;\bfL_R^{\langle - \infty \rangle}\bigr)\right).
$$
Using Lemma~\ref{lem:useful_G-homeo},
\begin{eqnarray*}
\lefteqn{H_{j}^{N_GC}\bigl(
\eub{N_GC} \times p_C^*W_GC \to  p_C^*W_GC;\bfL_R^{\langle - \infty \rangle}\bigr)}
& & 
\\
&\cong &
H_{j}^{N_GC}\bigl(
N_GC \times_C \res_{N_GC}^C\eub{N_GC} 
\to N_GC \times_C \pt;\bfL_R^{\langle - \infty \rangle}\bigr)
\\
&\cong &
H_{j}^{C}\bigl(
\res_{N_GC}^C\eub{N_GC} \to \pt;\bfL_R^{\langle - \infty \rangle}\bigr)
\\
&\cong &
H_{j}^{C}\bigl(EC \to \pt;\bfL_R^{\langle - \infty \rangle}\bigr).
\end{eqnarray*}
But this is zero by Remark~\ref{rem:role_of_type_I_and_II}. Hence, $E^2_{i,j} = 0$ 
for all $i,j \in \IZ$. This implies that for every free $W_GC$-$CW$-complex
$Y$, $H_*^{N_GC}\bigl(\eub{N_GC} \times p_C^*Y \to p_C^*Y;\bfL_R^{\langle -
\infty \rangle}\bigr)$ vanishes. In particular,
\[
H_n^{N_GC}\bigl(\eub{N_GC} \times p_C^* EW_GC 
\to  p_C^*EW_GC;\bfL_R^{\langle - \infty \rangle}\bigr)
=0\]
and
\[
H_n^{N_GC}\bigl(\eub{N_GC} \times p_C^*(W_GC \times_H EH) 
\to  p_C^*(W_GC \times_H EH);\bfL_R^{\langle - \infty \rangle}\bigr) =0.
\]
Therefore the Mayer-Vietoris sequence associated to the 
$W_GC$-pushout~\eqref{pushout_for_eub(W_GC)} yields an isomorphism
\begin{multline*}
\bigoplus_{H \in J_C} H_n^{N_GC}\bigl(\eub{N_GC} \times p_C^*W_GC/H 
\to  p_C^*W_GC/H ;\bfL_R^{\langle - \infty \rangle}\bigr)
\\
\xrightarrow{\cong}
H_*^{N_GC}\bigl(\eub{N_GC} \times p_C^*\eub{W_GC} 
\to  p_C^*\eub{W_GC};\bfL_R^{\langle - \infty \rangle}\bigr).
\end{multline*}
Since $p_C^{-1}(H) \cong D_{\infty}$, Lemma~\ref{lem:useful_G-homeo} implies:
\begin{eqnarray*}
\lefteqn{H_n^{N_GC}\bigl(\eub{N_GC} \times p_C^*W_GC/H 
\to  p_C^*W_GC/H ;\bfL_R^{\langle - \infty \rangle}\bigr)} \hspace{6mm}
& & 
\\
& \cong &
H_n^{N_GC}\bigl(N_GC \times_{p_C^{-1}(H)} \res_{N_GC}^{p_C^{-1}(H)} \eub{N_GC}  
\to N_GC \times_{p_C^{-1}(H)} \pt;\bfL_R^{\langle - \infty \rangle}\bigr)
\\
& \cong &
H_n^{p_C^{-1}(H)}\bigl(\res_{N_GC}^{p_C^{-1}(H)} \eub{N_GC}  
\to  \pt;\bfL_R^{\langle - \infty \rangle}\bigr)
\\
& \cong &
H_n^{p_C^{-1}(H)}\bigl(\eub{p_C^{-1}(H)}  
\to  \pt;\bfL_R^{\langle - \infty \rangle}\bigr)
\\
& \cong &
H_n^{D_{\infty}}\bigl(\eub{D_{\infty}}  
\to  \pt;\bfL_R^{\langle - \infty \rangle}\bigr)
\\
& = &
\UNil^{\langle - \infty \rangle}_n(D_{\infty};R).
\end{eqnarray*}
The projection $\eub{N_GC} \times p_C^*\eub{W_GC} \to \eub{N_GC}$
is a $N_GC$-homotopy equivalence, thus
\begin{eqnarray*}
\hspace{-8mm} \bigoplus_{H \in J_C} \UNil^{\langle - \infty \rangle}_n(D_{\infty};R)
&\xrightarrow{\cong} &
H_n^{N_GC}\bigl(\eub{N_GC} \to  p_C^*\eub{W_GC};\bfL_R^{\langle - \infty \rangle}\bigr).
\end{eqnarray*}
Therefore,
\begin{eqnarray*}
\bigoplus_{C \in I_2} \bigoplus_{H \in J_C} \UNil^{\langle - \infty \rangle}_n(D_{\infty};R)
& \cong &
H_n^G\bigl(\eub{G} \to \edub{G};\bfL^{\langle - \infty}_R\bigr),
\end{eqnarray*}
which implies that
$$
\left(\bigoplus_{F \in J} \cals_n^{\per,\langle - \infty \rangle}(F;R)\right)
\oplus
\left(\bigoplus_{C \in I_2} \bigoplus_{H \in J_C} 
\UNil^{\langle - \infty \rangle}_n(D_{\infty};R)\right)
\xrightarrow{\cong} 
\cals_n^{\per,\langle -\infty \rangle}(G;R).
$$
Theorem~\ref{the:L-theory_of_G_with_Aq_is_0} now follows from the fact that
$I_2$ is empty if $Q$ has odd order.
\end{proof}


\subsubsection{$K$-theory in the case $Q = \IZ/p$ for a prime $p$ and regular $R$}
\label{subsubsec:proof_K-theory_in_the_case_Q_is_Z/p_for_a_prime_p}

\begin{lemma}
\label{lem:case_pr-1(H)-torsionfree}
Let $f \colon G \to G'$ be a group homomorphism, $\calf$ be a family of
subgroups of $G$ and $R$ be a ring. Let $X$ be a $G'$-$CW$-complex such that
for every isotropy group $H' \subseteq G'$ of $X$ and every $n\in \IZ$,
$$H^{f^{-1}(H')}_n\left(E_{\calf\cap f^{-1}(H')}\bigl(f^{-1}(H')\to \pt\bigr);\bfK_R\right) = 0,$$
where $\calf\cap f^{-1}(H'):=\{ K\cap f^{-1}(H')\;|\;K\in\calf \}$. Then, for
every $n\in \IZ$,
$$H^G_n\bigl(E_\calf G \times f^*X
\to f^*X;\bfK_R\bigr) = 0.$$
\end{lemma}

\begin{proof}
Sending $X$ to $H^G_n\bigl(E_\calf G \times f^*X \to f^*X;\bfK_R\bigr)$
defines a $G'$-homology theory. There is an equivariant version of the
Atiyah-Hirzebruch spectral sequence converging to $H^G_{i+j}\bigl(E_\calf G
\times f^*X \to f^*X;\bfK_R\bigr)$ (see
Davis-L\"uck~\cite[Theorem~4.7]{Davis-Lueck(1998)}).  Its $E^2$-term is the
Bredon homology of $X$
$$E^2_{i,j} = H^{\Or(G')}_i(X;V_j),$$
where the coefficients are given by the covariant functor
$$V_j \colon \Or(G') \to \IZ\text{-}\MODULES$$
defined by
$$G'/H'
\mapsto 
H^G_j\bigl(E_\calf G \times f^*(G'/H') 
\to f^*(G'/H');
\bfK_R\bigr).$$
It suffices to show that the $E^2$-term is trivial for all $i,j$.
We will do this by showing that for every subgroup 
$H' \subseteq G'$ and every $j \in \IZ$,
$$H^G_j\bigl(E_\calf G \times f^*(G'/H') 
\to f^*(G'/H');
\bfK_R\bigr) = 0.$$
Using the induction structure and Lemma~\ref{lem:useful_G-homeo}
\begin{eqnarray*}
\lefteqn{H^G_j\bigl(E_\calf G \times f^*(G'/H') 
\to f^*(G'/H');
\bfK_R\bigr)}
& & 
\\ 
& \cong &
H^G_j\left(G \times_{f^{-1}(H')}  
\bigl(\res_G^{f^{-1}(H')}E_\calf G  \to \pt\bigr);\bfK_R\right) 
\\
& \cong &
H^G_j\left(G \times_{f^{-1}(H')}  
\bigl(E_{\calf\cap f^{-1}(H')}\bigl(f^{-1}(H')\bigr)  \to \pt\bigr);\bfK_R\right) 
\\
& \cong &
H^{f^{-1}(H')}_j\left(E_{\calf\cap f^{-1}(H')}\bigl(f^{-1}(H')\to \pt \bigr);\bfK_R\right),
\end{eqnarray*}
which is zero by assumption.
\end{proof}

The following general lemma will also be needed.

\begin{lemma} \label{lem;G_1_times_G_2_and_K}
Let $G_1$ and $G_2$ be two groups. Let $X$ be a $G_1$-$CW$-complex, and
let $\pr \colon G_1 \times G_2 \to G_1$ be projection.
Then, for every $n \in \IZ$, there are natural isomorphisms
\begin{eqnarray*}
H^{G_1 \times G_2}_n\bigl(\pr^*X;\bfK_R\bigr)
& \cong &
H^{G_1}_n\bigl(X;\bfK_{R[G_2]}\bigr);
\\
H^{G_1 \times G_2}_n\bigl(\pr^*X;\bfL^{\langle -\infty\rangle}_R\bigr)
& \cong &
H^{G_1}_n\bigl(X;\bfL^{\langle -\infty\rangle}_{R[G_2]}\bigr).
\end{eqnarray*}
\end{lemma}
\begin{proof}
We use the notation from \cite[Section~6]{Lueck-Reich(2005)}).  Let $\pr
\colon \Or(G_1 \times G_2) \to \Or(G_2)$ be the functor given by induction
with $\pr$. It sends $(G_1 \times G_2)/H$ to $G_1/\pr(H)$. Because of the
adjunction between induction and restriction
(see~\cite[Lemma~1.9]{Davis-Lueck(1998)})
and~\cite[Lemma~4.6]{Davis-Lueck(1998)}), it suffices to construct a weak
equivalence of covariant $\Or(G_1)$-spectra
$$
\pr_*\bigl(\bfK_R \circ \calg^{G_1 \times G_2}\bigr) \xrightarrow{\cong}
\bfK_{RG_2} \circ \calg^{G_1}.
$$
By definition, $\pr_*\bigl(\bfK_R \circ \calg^{G_1 \times G_2}\bigr)$ sends
$G_1/H$ to $\bfK_R\bigl(\calg^{G_1}(G_1) \times \widehat{G_2}\bigr)$, where
$\widehat{G_2}$ is the groupoid with one object and $G_2$ as its automorphism
group. The desired weak equivalence of covariant $\Or(G_1)$-spectra is now
obtained by comparing the definition of the category of free
$R\bigl(\calg^{G_1}(G_1) \times \widehat{G_2}\bigr)$-modules with the definition
of the category of free $R[G_2]\calg^{G_1}(G_1)$-modules.
\end{proof}

\begin{proof}[Proof of Theorem~\ref{the:K-theory_RG_Q_is_Z/p}]

By Theorem~\ref{the:FJC_hyperbolic_virtually_Zd} and~(\ref{relative_K_splitting}),
\begin{eqnarray*}
\Wh_n(G;R) 
\; \cong \;
H_n^G\bigl(EG\to \eub{G};\bfK_R\bigr) \oplus H_n^G\bigl(\eub{G}\to \edub{G};\bfK_R\bigr). 
\end{eqnarray*}
We begin by analyzing $H_n^G\bigl(EG\to \eub{G};\bfK_R\bigr)$.

Note that every non-trivial finite subgroup $H$ of $G$ is isomorphic to $\IZ/p$,
and is maximal. Furthermore, the normalizer $N_GH$ is isomorphic to $A^{\IZ/p}
\times H$. This can be seen as follows. The homomorphism $q \colon G \to Q = \IZ/p$ 
induces an injection $H \to \IZ/p$, which must be an isomorphism for
non-trivial $H$. Clearly $A^{\IZ/p}$ and $H$ belong to $N_GH$, and the subgroup
generated by $A^{\IZ/p}$ and $H$ is isomorphic to $A^{\IZ/p} \times H$. Thus, it
remains to show that $N_GH \subseteq A^{\IZ/p} \times H$. Let $t \in H$ be a
generator of $H$. Then every element in $G$ is of the form $at^i$ for some $i
\in \{0,1,2, \ldots p-1\}$ and some $a \in A$.  If $at^i \in N_GH$, then
$at^it(at^i)^{-1} = ata^{-1} = t^j$ for some $j \in \{0,1,2, \ldots p-1\}$.
Since $q(t) = q(ata^{-1}) = q(t^j) = q(t)^j$, $j = 1$. Thus $ata^{-1} = t$ and
$t^{-1}at = a$. This implies that $a \in A^{\IZ/p}$, and hence, $at^i \in
A^{\IZ/p} \times H$.

From~\cite[Corollary~2.10]{Lueck-Weiermann(2012)}, there is a $G$-pushout
$$\xymatrix{
     \coprod_{H\in J} G\times_{A^{\IZ/p} \times H}EA^{\IZ/p} \times EH
     \ar[d]^{\coprod_{H\in J} p} \ar[r]^-i & EG \ar[d] \\
     \coprod_{H\in J} G\times_{A^{\IZ/p} \times H}EA^{\IZ/p}\ar[r] & \eub{G}
 }
$$
which induces an isomorphism
\begin{multline*}
\bigoplus_{H\in J} H^G_n\bigl(G\times_{A^{\IZ/p} \times H}(EA^{\IZ/p} \times EH) 
\to G\times_{A^{\IZ/p} \times H}EA^{\IZ/p};\bfK_R\bigr) 
\\
\xrightarrow{\cong}
H^G_n\bigl(EG \to \eub{G};\bfK_R\bigr),
\end{multline*}
where $J$ is a complete system of representatives of the conjugacy classes of maximal finite subgroups of $G$.
Since $A^{\IZ/p} \cong \IZ^e$, the induction structure implies
\begin{eqnarray*}
\lefteqn{H^G_n\bigl(G\times_{A^{\IZ/p} \times H}(EA^{\IZ/p} \times EH) 
\to G\times_{A^{\IZ/p} \times H}EA^{\IZ/p};\bfK_R\bigr)}
\\
& \cong &
H^{A^{\IZ/p} \times H}\bigl(EA^{\IZ/p} \times EH 
\to EA^{\IZ/p};\bfK_R\bigr) 
\\
& \cong &
H^{H}\bigl(BA^{\IZ/p} \times (EH \to \pt);\bfK_R\bigr)
\\
& \cong &
\bigoplus_{i=0}^e H^{H}_{n-i}\bigl(EH \to \pt;\bfK_R\bigr)^{\binom{e}{i}}
\\
& \cong &
\bigoplus_{i=0}^e \Wh_{n-i}(H;R)^{\binom{e}{i}}.
\end{eqnarray*} 
Therefore,
\begin{eqnarray}
\bigoplus_{H\in J} \bigoplus_{i=0}^e \Wh_{n-i}(H;R)^{\binom{e}{i}}
& \xrightarrow{\cong} &
H^G_n\bigl(EG \to \eub{G};\bfK_R\bigr).
\label{HG_n(EG_to_eubG;bfK_R)_Q_is_Z/p}
\end{eqnarray}

Now we turn our attention to $H_n^G\bigl(\eub{G}\to \edub{G};\bfK_R\bigr)$. Let
$\overline{A} := A/A^{\IZ/p}$ and $\overline{G}:= G/A^{\IZ/p}$.  Then the exact
sequence~\eqref{1_to_A_to_G_to_Q_to_1} induces the exact sequence
\begin{eqnarray}
& 1 \to \overline{A} \to \overline{G} \xrightarrow{q} Q \to 1.
\label{1_to_overlineA_to_overlineG_to_Q_to_1}
\end{eqnarray}
Let $\calf$ be the family of subgroups $K$ of $G$ for which $\pr(K)$ is finite,
where $\pr:G\to \overline{G}$ is the quotient homomorphism. Let $\calvcyc$ be
the family of virtually cyclic subgroups of $G$, $\calvcyc_I$ be the family of
virtually cyclic subgroups of type I, and $\calf_1=\calf\cap\calvcyc_I$. The
Farrell-Jones Conjecture in algebraic $K$-theory is true for any group appearing
in $\calf$, since every element in $\calf$ is virtually finitely generated
abelian (see Theorem~\ref{the:FJC_hyperbolic_virtually_Zd}).  It is
straightforward to check that $\pr^*\eub{\overline{G}}$ is a model for
$\EGF{G}{\calf}$. By the Transitivity Principle
(see~\cite[Theorem~A.10]{Farrell-Jones(1993a)} and \cite[Theorem~65]{Lueck-Reich(2005)}) and
Remark~\ref{rem:role_of_type_I_and_II},
\begin{eqnarray}
H_n^G\bigl(\EGF{G}{\calf_1};\bfK_R\bigr) 
& \xrightarrow{\cong} &
H_n^G\bigl(\pr^*\eub{\overline{G}};\bfK_R\bigr).
\label{transitivity_F_1_to_F_0}
\end{eqnarray}

Since every virtually cyclic subgroup of type I in $\overline{G}$
is infinite cyclic, every element $K \in \calvcyc_I$ belongs to $\calf_1$, or
is infinite cyclic and $\{K' \mid K' \subseteq K \text{ and } K \in \calf_1\}$
consists of just the trivial group. Since $R$ is assumed to be regular, 
the map $H_n^{\IZ}(E\IZ;\bfK_R) \to H_n^{\IZ}(\pt;\bfK_R)$ is
bijective for every $n \in \IZ$. Therefore, the Transitivity Principle implies that
\begin{eqnarray}
H_n^G\bigl(\EGF{G}{\calf_1};\bfK_R\bigr) 
&
\cong
&
H_n^G\bigl(\EGF{G}{\calvcyc_I};\bfK_R\bigr).
\label{transitivity_F_1_to_vcyc_I}
\end{eqnarray}
By Remark~\ref{rem:role_of_type_I_and_II},
\begin{eqnarray}
H_n^G\bigl(\EGF{G}{\calvcyc_I};\bfK_R\bigr) 
&\xrightarrow{\cong} &
H_n^G\bigl(\edub{G};\bfK_R\bigr)
\label{transitivity_vcyc_I_to_vcyc}
\end{eqnarray}
is also a bijection. Thus,~\eqref{transitivity_F_1_to_F_0},
\eqref{transitivity_F_1_to_vcyc_I},~\eqref{transitivity_vcyc_I_to_vcyc} 
and the Five-Lemma imply that, for every $n\in \IZ$,
\begin{eqnarray}
H^G_n\bigl(\eub{G}\to \edub{G};\bfK_R\bigr)
\; \cong \;
H^G_n\bigl(\eub{G}\to \pr^*\eub{\overline{G}};\bfK_R\bigr).
\label{K_Z_p_second_piece}
\end{eqnarray}

Next we show that the conjugation action of $\IZ/p$ on $\overline{A}$
is free away from $0 \in \overline{A}$. Let $t$ be an element of $G$ that is
mapped by $q \colon G \to \IZ/p$ to a generator of $\IZ/p$. Let $a
\in A$ be given such that $\pr(t)\pr(a)\pr(t)^{-1}= \pr(a)$. 
Thus, there is a $b \in A^{\IZ/p}$ such that $tat^{-1} = ab$. 
If $\rho \colon A \to A$ denotes the conjugation action
with $t$, then $\rho(a) - a = b$ in $A$, where the
group operation in $A$ is written additively. Since
$$
p \cdot b 
= \sum_{i=0}^{p-1} \rho^i(b)
=  \sum_{i=0}^{p-1} \rho^i\bigl(\rho(a) - a\bigr)
= \sum_{i=0}^{p-1} \rho^{i+1}(a) - \rho^i(a)
=  \rho^p(a) - a
=  0,
$$
it follows that $\rho(a) - a=b=0$. This implies that $a \in A^{\IZ/p}$, and thus
$\overline{a} = 0$ in $\overline{A}$. Therefore, the conjugation action of
$\IZ/p$ on $\overline{A}$ is free away from $0 \in \overline{A}$.

Let $\overline{J}$ be a complete system of representatives of 
the conjugacy classes of maximal finite subgroups of $\overline {G}$. 
From~\cite[Lemma~6.3]{Lueck-Stamm(2000)} 
and~\cite[Corollary~2.11]{Lueck-Weiermann(2012)}, there is a $\overline{G}$-pushout
\begin{eqnarray}
& \xymatrix{
     \coprod_{\overline{H}\in \overline{J}} \overline{G}\times_{\overline{H}}E\overline{H}
     \ar[d] \ar[r]^-i & E\overline{G} \ar[d] \\
     \coprod_{\overline{H} \in \overline{J}} \overline{G}\times_{\overline{H}} \pt\ar[r] & \eub{\overline{G}}
 }
&
\label{overlineG_pushout}
\end{eqnarray}

Sending a $\overline{G}$-$CW$-complex $\overline{X}$ to the $\IZ$-graded abelian
group $H^G_*\bigl(\eub{G} \times \pr^*\overline{X} \to
\pr^*\overline{X};\bfK_R\bigr)$ defines a $\overline{G}$-homology theory. Since
$\pr^{-1}(\{1\}) = A^{\IZ/p}\cong \IZ^e$ and $R$ is regular,
Theorem~\ref{the:Zd} and Lemma~\ref{lem:case_pr-1(H)-torsionfree} imply that if
$\overline{X}$ is a free $\overline{G}$-$CW$-complex, then 
$H^G_*\bigl(\eub{G} \times \pr^*\overline{X} \to \pr^*\overline{X};\bfK_R\bigr)=0$. Therefore,
because $ \coprod_{\overline{H}\in \overline{J}}
\overline{G}\times_{\overline{H}}E\overline{H}$ and $E\overline{G}$ are free
$\overline{G}$-$CW$-complexes, the Mayer-Vietoris sequence associated to the
$\overline{G}$-pushout~\eqref{overlineG_pushout} yields an isomorphism
\begin{multline}
\bigoplus_{\overline{H} \in \overline{J}} 
H^G_n\bigl(\eub{G} \times \pr^*\overline{G}/\overline{H} 
\to \pr^*\overline{G}/\overline{H};\bfK_R\bigr)
\\
\cong 
H^G_n\bigl(\eub{G} \times \pr^*\eub{\overline{G}} \to \pr^*\eub{\overline{G}};\bfK_R\bigr).
\label{HG_n(eubGto_euboverlineG;K_R)_(0)}
\end{multline}

Notice that if $\pr^{-1}(\overline{H})$ is torsion-free, then it is isomorphic
to $\IZ^e$. If $\overline{H}$ is trivial, then this follows from
$\pr^{-1}(\{1\}) = A^{\IZ/p}$. If $\overline{H}$ is a non-trivial finite
subgroup of $\overline{G}$, choose an element $t \in \pr^{-1}(\overline{H})$
such that $\pr(t)$ is a generator of $\overline{H}$. Then every element in
$\pr^{-1}(\overline{H})$ is of the form $at^u$ for $a \in A^{\IZ/p}$ and $u \in
\{0,1,2,\ldots , p-1\}$. For two such elements $at^u$ and $bt^v$
$$at^ubt^v = abt^ut^v = bat^vt^u = bt^vat^u.$$
Hence $\pr^{-1}(\overline{H})$ is a torsion-free abelian group containing
$\IZ^e$ as subgroup of finite index, and so $\pr^{-1}(\overline{H}) \cong
\IZ^e$. Thus, if $\pr^{-1}(\overline{H})$ is torsion-free, then
Theorem~\ref{the:Zd} and Lemma~\ref{lem:case_pr-1(H)-torsionfree} imply that
\[
H^G_n\bigl(\eub{G} \times \pr^*\overline{G}/\overline{H} \to
\pr^*\overline{G}/\overline{H}; \bfK_R\bigr) =0.
\]
Therefore, if $\overline{J}'$ is the subset of $\overline{J}$ consisting of
those elements $\overline{H} \in J$ for which $\pr^{-1}(\overline{H})$ is not
torsion-free, then~\eqref{HG_n(eubGto_euboverlineG;K_R)_(0)} becomes
\begin{multline}
\bigoplus_{\overline{H} \in \overline{J}'} H^G_n\bigl(\eub{G} \times
\pr^*\overline{G}/\overline{H} \to \pr^*\overline{G}/\overline{H};\bfK_R\bigr)
\\
\cong H^G_n\bigl(\eub{G} \times \pr^*\eub{\overline{G}} \to
\pr^*\eub{\overline{G}};\bfK_R\bigr).
\label{HG_n(eubGto_euboverlineG;K_R)_(1)}
\end{multline}

Let $\overline{H} \in \overline{J}'$. Then $\pr^{-1}(\overline{H}) \cong
A^{\IZ/p} \times \overline{H}$. Since $A^{\IZ/p}\cong \IZ^e$, the induction
structure, Lemma~\ref{lem:useful_G-homeo}, Lemma~\ref{lem;G_1_times_G_2_and_K}
and Theorem~\ref{the:Zd} imply
\begin{eqnarray*}
\lefteqn{H^G_n\bigl(\eub{G} \times \pr^*\overline{G}/\overline{H} 
  \to \pr^*\overline{G}/\overline{H};\bfK_R\bigr)}
\\
& \cong &
H^{G}_n
\biggl(G \times_{A^{\IZ/p} \times \overline{H}} \bigl(\res_G^{A^{\IZ/p} \times \overline{H}} \eub{G} \to \pt\bigr);
\bfK_R\biggr)
\\
& \cong &
H^{A^{\IZ/p} \times \overline{H}}_n
\bigl(\res_G^{A^{\IZ/p} \times \overline{H}} \eub{G} \to \pt;\bfK_R\bigr)
\\
& \cong &
H^{A^{\IZ/p} \times \overline{H}}_n
\bigl(EA^{\IZ/p} \to \pt;\bfK_R\bigr)
\\
& \cong &
H^{A^{\IZ/p}}_n \bigl(EA^{\IZ/p} \to \pt;\bfK_{R[\overline{H}]}\bigr)
\\
& \cong &
\Wh_n(A^{\IZ/p};R[\IZ/p])
\\
& \cong &
\bigoplus_{C \in  \calmicyc(A^{\IZ/p})} \;\bigoplus_{i = 0}^{e-1}  
\bigl(\NK_{n-i}(R[\IZ/p]) \oplus \NK_{n-i}(R[\IZ/p])\bigr)^{\binom{e-1}{i}}.
\end{eqnarray*}
Since the projection $\eub{G}\times \pr^*\eub{\overline{G}} \to \eub{G}$ is a
$G$-homotopy equivalence,~\eqref{HG_n(eubGto_euboverlineG;K_R)_(1)} implies that
\begin{multline*}
\bigoplus_{\overline{H} \in \overline{J}'} \bigoplus_{C \in
  \calmicyc(A^{\IZ/p})} \;\bigoplus_{i = 0}^{e-1} \bigl(\NK_{n-i}(R[\IZ/p])
\oplus \NK_{n-i}(R[\IZ/p])\bigr)^{\binom{e-1}{i}}
\\
\cong H^G_n\bigl(\eub{G}\to \pr^*\eub{\overline{G}};\bfK_R\bigr).
\end{multline*}
Together with~\eqref{K_Z_p_second_piece}, this produces the isomorphism
\begin{multline*}
\bigoplus_{\overline{H} \in \overline{J}'} \bigoplus_{C \in
  \calmicyc(A^{\IZ/p})} \; \left(\bigoplus_{i = 0}^{e-1}
  \bigl(\NK_{n-i}(R[\IZ/p]) \oplus
  \NK_{n-i}(R[\IZ/p])\bigr)^{\binom{e}{i}}\right) \\ \cong
H^G_n\bigl(\eub{G}\to \edub{G};\bfK_R\bigr).
\end{multline*}

To complete the proof of the theorem, we must show that there is a bijection
between the sets $J$ and $\overline{J}'$. Send an element $H \in J$ to the
element in $H' \in \overline{J}'$ that is uniquely determined by the property
that $\pr(H)$ and $H'$ are conjugated in $\overline{G}$.  This map is
well-defined and injective because two non-trivial finite subgroups $H_1$ and
$H_2$ of $G$ are conjugate if and only if $\pr(H_1)$ and $\pr(H_2)$ are
conjugate in $\overline{G}$. The map is surjective, since for any element
$\overline{H} \in \ \overline{J}'$ there exists a finite subgroup $H \subseteq
G$ with $\overline{H} = \pr(H)$. This finishes the proof of
Theorem~\ref{the:K-theory_RG_Q_is_Z/p}.
\end{proof}


\subsubsection{$L$-theory in the case $Q = \IZ/p$ for an odd prime $p$}
\label{subsubsec:proof_L-theory_in_the_case_Q_is_Z/p_for_an_odd_prime_p}

\begin{proof}[Proof of Theorem~\ref{the:L-theory_RG_Q_is_Z/p}]
An argument completely analogous to the one that
established~\eqref{HG_n(EG_to_eubG;bfK_R)_Q_is_Z/p} in the proof of
Theorem~\ref{the:K-theory_RG_Q_is_Z/p} shows that
\begin{eqnarray}
  \bigoplus_{H\in J} \bigoplus_{i=0}^e \cals^{\per,\langle -\infty \rangle}_{n-i}(H;R)^{\binom{e}{i}}
  & \xrightarrow{\cong} &
  H^G_n\bigl(EG \to \eub{G};\bfL_R^{\langle -\infty \rangle}\bigr).
  \label{HG_n(EG_to_eubG;bfL_R)_Q_is_Z/p}
\end{eqnarray}
Since $p$ is odd, every infinite virtually cyclic subgroup of $G$ is of type
I. Thus the result follows from Theorem~\ref{the:FJC_hyperbolic_virtually_Zd},
Remark~\ref{rem:role_of_type_I_and_II}
and~\eqref{HG_n(EG_to_eubG;bfL_R)_Q_is_Z/p}.
\end{proof}

\begin{proof}[Proof of Theorem~\ref{the:L-theory_ZG_Q_is_Z/p}]
There are natural maps $\bfL_{\IZ}^{\langle j \rangle} \to \bfL_{\IZ}^{\langle
  j+1 \rangle}$ of covariant functors from $\GROUPOIDS$ to $\SPECTRA$ for $j
\in \{2,1,0,-1, \ldots\}$.  Let $\bfL_{\IZ}^{\langle j+1,j\rangle}$ be the
covariant functor from $\GROUPOIDS$ to $\SPECTRA$ that assigns to a groupoid
$\calg$ the homotopy cofiber of the map of spectra obtained by evaluating
$\bfL_{\IZ}^{\langle j \rangle} \to \bfL_{\IZ}^{\langle j+1 \rangle}$ at
$\calg$. We obtain equivariant homology theories
$H_*^?\bigl(-;\bfL_{\IZ}^{\langle j+1 \rangle}\bigr)$,
$H_*^?\bigl(-;\bfL_{\IZ}^{\langle j \rangle}\bigr)$ and
$H_*^?\bigl(-;\bfL_{\IZ}^{\langle j+1,j \rangle}\bigr)$ (see
L\"uck-Reich~\cite[Section~6]{Lueck-Reich(2005)}), such that for any group $G$
and any $G$-$CW$-complex $X$, there is a long exact sequence
\begin{multline}
 \cdots \to H_n^G\bigl(X;\bfL_{\IZ}^{\langle j+1 \rangle}\bigr) \to
 H_n^G\bigl(X;\bfL_{\IZ}^{\langle j \rangle}\bigr) \to 
 H_n^G\bigl(X;\bfL_{\IZ}^{\langle j+1,j \rangle}\bigr) 
 \\ \to H_{n-1}^G\bigl(X;\bfL_{\IZ}^{\langle j+1 \rangle}\bigr) \to
 H_{n-1}^G\bigl(X;\bfL_{\IZ}^{\langle j \rangle}\bigr) \to \cdots.
 \label{pseudo_Rothenberg-sequence}
\end{multline}
If $X = \pt$, then 
\begin{eqnarray*}
H_n^G\bigl(\pt;\bfL_{\IZ}^{\langle j+1 \rangle}\bigr) 
& \cong & 
L_n^{\langle j+1 \rangle}(\IZ G);
\\
H_n^G\bigl(\pt;\bfL_{\IZ}^{\langle j \rangle}\bigr) 
& \cong & 
L_n^{\langle j \rangle}(\IZ G);
\\
H_n^G\bigl(\pt;\bfL_{\IZ}^{\langle j+1,j \rangle}\bigr) 
& \cong & 
\widehat{H}^n(\IZ/ 2;\widetilde{K}_{j}(\IZ G)),
\end{eqnarray*}
and the sequence~\eqref{pseudo_Rothenberg-sequence}
can be identified with the Rothenberg sequence~\eqref{Rothenberg_sequence}.
Since $L_n^{\langle j+1 \rangle }(\IZ)  \xrightarrow{\cong} L_n^{\langle j \rangle}(\IZ)$ is 
a bijection for every $n \in \IZ$ by~\eqref{Rothenberg_sequence}, 
a spectral sequence argument shows that the map
$$H_n\bigl(BG;\bfL^{\langle j+1 \rangle}(\IZ)\bigr)
\xrightarrow{\cong} 
H_n\bigl(BG;\bfL^{\langle j \rangle}(\IZ)\bigr)
$$
is a bijection for every $n \in \IZ$. Thus,
$H_n^G\bigl(EG;\bfL^{\langle j+1,j \rangle}_{\IZ}\bigr) 
\cong H_n\bigl(BG;\bfL^{\langle j+1,j \rangle}(\IZ)\bigr)$ is zero for all $n \in \IZ$.
Hence, for every $n \in \IZ$, there is an isomorphism
$$H_n^G\bigl(EG \to \pt;\bfL^{\langle j+1,j \rangle}_{\IZ}\bigr) 
\xrightarrow{\cong} 
H_n^G\bigl(\pt;\bfL^{\langle j+1,j \rangle}_{\IZ}\bigr).$$
This implies that there is also a Rothenberg sequence for the structure sets
\begin{multline}
\cdots \to \cals_n^{\per,\langle j+1 \rangle}(G;\IZ) 
\to \cals_n^{\per,\langle j \rangle}(G;\IZ) \to
\widehat{H}^n(\IZ/ 2;\widetilde{K}_{j}(\IZ G)) 
\\ 
\to   \cals_{n-1}^{\per,\langle j+1 \rangle}(G;\IZ) 
\to \cals_{n-1}^{\per,\langle j \rangle}(G;\IZ) \to \cdots.
\label{Rothenberg_sequence_for_cals}
\end{multline}

Since $\NK_n(\IZ[\IZ/p]) = 0$ for $n \le 1$ 
(see Bass-Murthy~\cite{Bass-Murthy(1967)}), Theorem~\ref{the:K-theory_RG_Q_is_Z/p} 
implies that, for every $n \le 1$,
\begin{eqnarray}
\bigoplus_{H\in J}  \bigoplus_{i=0}^e  \Wh_{n-i}(H;\IZ)^{\binom{e}{i}}
& \xrightarrow{\cong} & 
\Wh_n(G;\IZ).
\label{iso_for_Wh_n(G)}
\end{eqnarray}
For every $H \in J$ and $j \le -2$, $K_j(\IZ H) = 0$ by Carter~\cite{Carter(1980)}.
Hence $K_j(\IZ G)$ vanishes for $j \le -2$ by~\eqref{iso_for_Wh_n(G)}. 
Theorem~\ref{the:L-theory_RG_Q_is_Z/p} and~\eqref{Rothenberg_sequence_for_cals} 
imply that the map
$$\bigoplus_{H\in J}  \bigoplus_{i=0}^e \cals_n^{\per,\langle j \rangle}(H;\IZ)^{\binom{e}{i}}
\to \cals_n^{\per,\langle j \rangle}(G;\IZ)$$
is bijective for all $n \in \IZ$ and $j \in \{-1,-2, \ldots\} \amalg \{-\infty\}$.

Next we use an inductive argument to show that this is also true for $j
=0,1$. By taking the direct sum of the Rothenberg
sequences~\eqref{Rothenberg_sequence_for_cals} for the various elements $H \in
J$ and mapping it to~\eqref{Rothenberg_sequence_for_cals} for $G$, one obtains
the following commutative diagram.
$$\xymatrix{
\vdots \ar[d] 
&
\vdots \ar[d] 
\\
\bigoplus_{H\in J}  \bigoplus_{i=0}^e 
\cals_{n-i}^{\per,\langle j+1 \rangle}(H;\IZ)^{\binom{e}{i}} \ar[r] \ar[d]
&
\cals_{n-i}^{\per,\langle j +1 \rangle}(G;\IZ) \ar[d]
\\
\bigoplus_{H\in J}  \bigoplus_{i=0}^e 
\cals_{n-i}^{\per,\langle j \rangle}(H;\IZ)^{\binom{e}{i}} \ar[r] \ar[d]
&
\cals_n^{\per,\langle j \rangle}(G;\IZ) \ar[d]
\\
\bigoplus_{H\in J}  \bigoplus_{i=0}^e 
\widehat{H}^{n-i}(\IZ/ 2;\widetilde{K}_{j}(\IZ H))^{\binom{e}{i}} 
\ar[r]^-{\cong}  \ar[d]
&
\widehat{H}^n(\IZ/ 2;\widetilde{K}_{j}(\IZ G)) \ar[d]
\\
\bigoplus_{H\in J}  \bigoplus_{i=0}^e 
\cals_{n-1-i}^{\per,\langle j+1 \rangle}(H;\IZ)^{\binom{e}{i}} \ar[r] \ar[d]
&
\cals_{n-1}^{\per,\langle j +1 \rangle}(G;\IZ) \ar[d]
\\
\bigoplus_{H\in J}  \bigoplus_{i=0}^e 
\cals_{n-1-i}^{\per,\langle j \rangle}(H;\IZ)^{\binom{e}{i}} \ar[r] \ar[d]
&
\cals_{n-1}^{\per,\langle j \rangle}(G;\IZ) \ar[d]
\\
\vdots
&
\vdots
}
$$
The middle horizontal arrow is an isomorphism since it is induced by the
isomorphism~\eqref{iso_for_Wh_n(G)}, which is compatible with the involutions.
The Five-Lemma implies that the maps
$$\bigoplus_{H\in J}  \bigoplus_{i=0}^e 
\cals_{n-i}^{\per,\langle j \rangle}(H;\IZ)^{\binom{e}{i}} \to
\cals_n^{\per,\langle j \rangle}(G;\IZ)$$ are bijective for all $n \in \IZ$ if
and only if the maps
$$\bigoplus_{H\in J}  \bigoplus_{i=0}^e \cals_{n-i}^{\per,\langle j+1 \rangle}(H;\IZ)^{\binom{e}{i}}
\to \cals_n^{\per,\langle j+1 \rangle}(G;\IZ)$$ are bijective for all $n \in
\IZ$. This takes care of the decorations $j=1$ and $j=0$.  A similar argument,
where $\widetilde{K}_{j}(\IZ H)$ and $\widetilde{K}_{j}(\IZ G)$ are replaced by $\Wh(H)$ 
and $\Wh(G)$ in the Rothenberg sequence, shows that this is also
true for the decoration $s$, since it is true for the decoration $h$ which is
the the same as $\langle 1 \rangle$.  This completes the proof of the theorem.
\end{proof}


\typeout{------------------------------------------ References ---------------------------------------}


\end{document}